\documentclass[11pt, a4paper]{amsart}

\usepackage[english]{babel}
\usepackage{amsmath}
\usepackage{amssymb}
\usepackage{hyperref} 
\hypersetup{
    colorlinks=true,
    linkcolor=blue,
    filecolor=magenta,
    urlcolor=cyan,
    pdftitle={Neural Learning of Fast Matrix Multiplication Algorithms: a Geometric Perspective},
    pdfpagemode=FullScreen,
}
\usepackage{amsthm}
\usepackage{latexsym}
\usepackage{mathtools}
\usepackage{thmtools}
\usepackage{amsfonts}
\usepackage{mathrsfs}
\usepackage{textcomp}
\usepackage{graphicx}
\usepackage{setspace}
\usepackage{nicefrac}
\usepackage{indentfirst}
\usepackage{enumerate}
\usepackage{wasysym}
\usepackage[normalem]{ulem}
\usepackage{upgreek}
\usepackage{paralist}
\usepackage{xcolor}
\usepackage{ulem}

\usepackage{etoolbox}
\usepackage{mathdots}
\usepackage{wrapfig}
\usepackage{floatflt}
\usepackage{tensor}
\usepackage{parskip}
\usepackage{lscape}
\usepackage{stmaryrd}
\usepackage[enableskew]{youngtab}
\usepackage{ytableau}
\usepackage{pifont}
\usepackage{braket}
\usepackage{derivative}
\usepackage[labelfont=bf, font={small}]{caption}[2005/07/16]
\usepackage{subcaption}

\usepackage[top=3cm, bottom=3cm, left=3cm, right=3cm]{geometry}
\usepackage[normalem]{ulem} % For \sout command
\usepackage{cleveref}

\usepackage{nicematrix} %This package is based on the package array. It creates PGF/TikZ nodes under the cells of the array and uses these nodes to provide functionalities to construct tabulars, arrays and matrices.

\usepackage{tikz}
\usetikzlibrary{cd}
\usetikzlibrary{arrows,automata}
\usetikzlibrary{matrix}
\usetikzlibrary{positioning}
\usetikzlibrary{snakes}
\usetikzlibrary{calc}

\definecolor{red}{rgb}{1,0,0}

\newcommand{\bbC}{\mathbb{C}}

\newcommand{\bbK}{\mathbb{K}}

\renewcommand{\phi}{\varphi}

\renewcommand{\bar}[1]{\overline{#1}}

\usepackage{algorithm}
\usepackage{algpseudocode}
\usepackage{booktabs}

\usepackage[all, cmtip]{xy}

\theoremstyle{definition}
\newtheorem{theorem}{Theorem}

\newtheorem{lemma}[theorem]{Lemma}

\newtheorem{definition}[theorem]{Definition}

\newtheorem{remark}[theorem]{Remark}

\newtheorem{procedure}[theorem]{Procedure}
\crefname{procedure}{Procedure}{Procedures}
\Crefname{procedure}{Procedure}{Procedures}

\title[]{Neural Learning of Fast Matrix Multiplication Algorithms: a Geometric Perspective}

\author[]{Paolo Andreini, Alessandra Bernardi, Monica Bianchini, Barbara Toniella Corradini, Sara Marziali, Giacomo Nunziati, Franco Scarselli}

\date{}

\address[P. Andreini, M. Bianchini, G. Nunziati, F. Scarselli]{Dipartimento di Ingegneria dell'Informazione e Scienze Matematiche, Universit\`a di Siena, Italy}
\address[A. Bernardi]{Dipartimento di Matematica, Universit\`a di Trento, Italy}
\address[B. T. Corradini]{Italian Institute of Technology, Genoa, Italy}
\address[S. Marziali]{Max Planck Institute of Molecular Cell Biology and Genetics (MPI-CBG); Center for Systems Biology Dresden (CSBD); Faculty of Mathematics, TU Dresden, Germany}

\email[(Andreini)]{paolo.andreini@unisi.it}
\email[(Bianchini)]{monica.bianchini@unisi.it}
\email[(Bernardi)]{alessandra.bernardi@unitn.it}
\email[(Corradini)]{barbaracorradini5@gmail.com}
\email[(Marziali)]{marziali@mpi-cbg.de}
\email[(Nunziati)]{giacomo.nunziati@student.unisi.it}
\email[(Scarselli)]{franco.scarselli@unisi.it}


\begin{document}

\begin{abstract}
We combine neural learning and geometric diagnostics to investigate the rank and border rank of the \(3\times3\) matrix-multiplication tensor \(M_3\). Using the Bhattacharya--Mesner product, we derive \textsc{MatMulNet}, a differentiable architecture encoding bilinear maps.
Constrained refinement yields a length-\(23\) numerical solution with machine-precision full-tensor residual and bounded decomposition weight. Along the explored length-\(22\) branch, the residual instead decreases as the gauge-invariant weight grows. This contrast supports \(R_{\mathbb R}(M_3)=23\) as a working hypothesis. A Laurent-parametrised extension probes border rank, giving a tentative, protocol-dependent motivation to investigate \(r=19\) and the adjacent \(r=18\) case.
For \(M_4\), the same diagnostics distinguish bounded numerical solutions recovered from perturbed known schemes, structured residual plateaus, and weight-growing
approximations, revealing the dependence of the search on initialisation.
Together, the experiments connect approximation accuracy with decomposition growth and identify concrete targets for exact or certified investigation.

\end{abstract}
\maketitle

\section{Introduction}

The search for fast matrix multiplication algorithms connects a concrete computational task with fundamental questions about tensor decompositions.
Strassen's $7$-multiplication algorithm for two \(2\times2\) matrices showed that algebraic identities at a small matrix size can change the asymptotic complexity of a basic operation~\cite{strassen1969gaussian}.
A bilinear algorithm with \(r\) scalar multiplications is equivalent to a length-\(r\) decomposition of the matrix-multiplication tensor; allowing limits of such decompositions leads to border rank. These two notions capture exact and approximate bilinear complexity, and both enter the characterisation of the matrix-multiplication exponent through block recursion and the theory of approximate algorithms~\cite{bcs1997,bini1980relations,schoenhage1981}.
At fixed matrix sizes, however, determining the minimum decomposition length remains a separate and remarkably delicate problem.

The first unresolved square case is \(M_3:=\langle3,3,3\rangle\). Over the fields relevant here, the known bounds for the rank $R(M_3)$ and the border-rank $\underline{R}(M_3)$ are respectively 
\[
19\leq R_{\mathbb C}(M_3)\leq R_{\mathbb R}(M_3)\leq23,
\qquad
17\leq\underline R_{\mathbb C}(M_3)
\leq\underline R_{\mathbb R}(M_3)\leq20.
\]
The rank bounds follow from Bl\"aser's lower bound and Laderman's 23-multiplication construction~\cite{blaser2003,laderman1976}; the border-rank bounds follow from the work of Conner--Harper--Landsberg and
Smirnov~\cite{conner-harper-landsberg-2023,smirnov}.
Here real border rank uses Euclidean closure, whereas complex border rank uses Zariski closure. Both the known length-23 decomposition and the known length-20 degeneration are available over the reals, the parameter field of our numerical experiments. In particular, whether a real length-$22$ decomposition exists remains unresolved.

Numerical searches offer a way to investigate this problem, but their interpretation requires more than a small reconstruction error. The image of the length-\(r\) decomposition map need not be closed: the represented tensors may converge while their rank-one summands become unbounded and cancel in the sum~\cite{deSilvaLim2008illposed}.
Consequently, an improving fit can reflect either approach to an exact finite decomposition or approach to a limiting configuration. Indeed, the known real border-rank upper bound already implies that the infimum of the unrestricted approximation error to \(M_3\) is zero for every \(r\geq20\). What a finite numerical experiment must distinguish is
therefore not just small from large error, but the different ways in which small error is attained.

Computational approaches to matrix multiplication explore several complementary representations of the problem. SAT-based searches have produced thousands of inequivalent length-23 schemes~\cite{heule2021newways}, while flip-graph methods search among exact schemes through local transformations~\cite{kauersMoosbauer2023flip}.
AlphaTensor uses reinforcement learning to search for decompositions in finite coefficient spaces~\cite{fawzi2022alphatensor}.
These methods demonstrate the effectiveness of automated discovery, but also the richness of the solution space: a learned length-23 representation need not resemble Laderman's original one.

Continuous optimisation has equally direct precedents. Tichavsk\'y, Phan and Cichocki used constrained Levenberg--Marquardt optimisation for matrix-multiplication tensors and studied reconstruction error as the permitted coefficient norm increases, including a comparison between lengths 22 and 23 for \(M_3\)~\cite{tichavsky2017numerical}.
Vermeylen and Van Barel developed an augmented-Lagrangian approach to sparse decompositions with bounded coefficients, obtaining decompositions and parameter families used to improve numerical stability~\cite{vermeylen2025stability}.
These works establish constrained numerical fitting as a useful tool for exploring short decompositions and the trade-off between approximation error and coefficient size.

The bilinear structure also admits direct learning from multiplication examples. Elser's multiplier network has the same three-stage bilinear structure---linear combinations of the two inputs, componentwise products, and linear output recombination. It recovered seven-term models for \(2\times2\) multiplication and 23-term models for \(3\times3\) multiplication, and exhibited slow convergence with growing weights at length 22~\cite{elser2016network}.
The StrassenNets architecture subsequently used trainable sum-product networks with a prescribed multiplication count and ternary  weights~\cite{tschannen2018strassennets}.
A more recent example combines backpropagation with ternarisation and symbolic verification to obtain an exact length-23 scheme~\cite{stapleton2025rank23}.
Numerical algebraic geometry provides another route from computation to decomposition questions: Bernardi, Daleo, Hauenstein and Mourrain developed homotopy-based methods for ranks, border ranks and decompositions using numerical descriptions of joins~\cite{BERNARDI201778}.

Starting from this background, our aim is to investigate what geometric information can be extracted from a learning-based search before an exact certificate is available. We combine an explicit algebraic description of the computation graph with experiments across decomposition lengths, full-tensor verification, and diagnostics of the size of the rank-one summands. The object of study is both the multiplication tensor and the behaviour of numerical searches near its possible short decompositions. The resulting evidence is intended to guide subsequent exact work, rather than to identify an algebraic invariant from an optimisation threshold.

We first construct tensor-network realisations of classical and Strassen \(2\times2\) multiplication using the Bhattacharya--Mesner product and the total-tensor formalism of directed networks~\cite{MESNER1990204,MESNER1994595,
chiantini2025producttensorsdescriptionnetworks}.
The blow and forget operations make explicit how the graph and its activation tensors encode the bilinear computation.
This construction leads to \textsc{MatMulNet}, whose real parameters \(\theta=(H,K,F)\) represent
\[
T_\theta=\sum_{s=1}^{r}h_s\otimes k_s\otimes f_s,
\qquad
H,K\in\mathbb R^{n^2\times r},\quad
F\in\mathbb R^{r\times n^2}.
\]
The three-factor parametrisation is classical; its role here is to link the network construction explicitly to a decomposition that can be inspected after training. Learning minimises the error on matrix-product examples. For our choice of random inputs, the expected squared error on matrix products is proportional to the squared Frobenius distance between the learned tensor and the matrix-multiplication tensor. Thus, minimising the expected loss is equivalent to minimising
the error on the tensor itself. 

The geometric analysis is organised around the pair
\[
\rho(\theta)=\frac{\|T_\theta-M_3\|_F}{\|M_3\|_F},
\qquad
W(\theta)=\sum_{s=1}^{r}
\|h_s\|_2\,\|k_s\|_2\,\|f_s\|_2.
\]
Unlike an individual coefficient norm, \(W\) is unchanged by the termwise rescalings that preserve each rank-one summand. With Euclidean norms it is the decomposition growth factor studied by Dai and Lim in their analysis of bilinear stability~\cite{daiLim2023stability}; it depends on the chosen decomposition, not only on the represented tensor.
We use it as a diagnostic of approximation at fixed decomposition length.
An elementary compactness argument, stated in \Cref{lem:bounded-weight}, shows that convergence with uniformly bounded
\(W\) forces the limit to have rank at most \(r\).
Conversely, approximating a tensor of higher rank with \(r\) terms requires \(W\) to diverge. This gives a precise mathematical basis for examining residual and weight together, while leaving open the possibility that a weight-growing numerical branch coexists with a different bounded representation of the same target.

The experiments reveal a distinction between accessibility of low-error fits and their subsequent behaviour under refinement.
The \(2\times2\) case calibrates the learning procedure against the known rank-seven regime without supplying Strassen's coefficients.
For \(M_3\), seven independent runs at each \(r\in\{16,\ldots,24\}\) produce no validation loss below \(10^{-3}\) at \(r\leq21\). 
Such fits occur in one of seven runs at \(r=22\), four at \(r=23\), and five at \(r=24\). We then evaluate selected candidates on all tensor coefficients and apply coefficient-constrained refinement.
Along the explored length-22 branch, the residual decreases approximately as \(W^{-1.02}\), reaching \(4.5\times10^{-6}\) at
\(W\simeq6.3\times10^5\). At length 23, refinement reaches a full-tensor relative residual of \(5.3\times10^{-16}\) with
\(W\simeq76\), giving a bounded numerical solution at double-precision resolution. The recorded refinement curves describe locally attained errors, not global minima under the coefficient constraints. 
A projected pseudo-arclength continuation of a length-22 candidate provides a complementary local exploration; discrepancies on the full tensor are checked separately from those on the projected equations.

This contrast develops the earlier observations of Elser and Tichavsk\'y--Phan--Cichocki through a decomposition-weight diagnostic and a joint reading of the sweep, refinement and continuation experiments.
It supports \(R_{\mathbb R}(M_3)=23\) as a working hypothesis and motivates further investigation of length 22. It does not certify non-existence of a length-22 decomposition. The information supplied by the experiments is more specific: they identify low-error candidates and quantify the cost, in summand size, of improving the candidates reached
by the search.

We also introduce an explicitly Laurent-parametrised extension in which all three factors depend on \(\varepsilon\). This allows the search to represent candidate border-rank families rather than only individual finite decompositions. The experiments cover \(r\in\{16,\ldots,21\}\) and three exponent ranges. The best validation losses decrease with \(r\), but remain above \(10^{-2}\), without a transition comparable to that of the finite-rank sweep.
The profile motivates only a tentative priority for \(r=19\), together with the neighbouring \(r=18\) case, in subsequent exact searches; this priority is specific to the tested protocol, not an estimate of the border rank. The distinction is essential because fitting at finite \(\varepsilon\) does not by itself enforce regularity or the desired limit at zero. For a fixed learned family, the relevant algebraic checks are cancellation of the negative-power coefficients and agreement of the constant coefficient with \(M_3\).

We also test the method at larger scale on \(M_4=\langle4,4,4\rangle\). The residual-weight analysis remains informative there: it separates bounded numerical solutions recovered from perturbed copies of a known length-\(48\) scheme, structured residual plateaus, and approximations whose weight grows as the residual decreases. No random start meets the numerical solution, which characterises the scaling of the search and shows the role of the distance of the initial point from a solution.

The resulting programme connects numerical exploration with concrete algebraic tasks. A bounded near-solution is a candidate for exact reconstruction or certified local analysis; a weight-growing branch motivates the search for a structured limiting family and the verification of its coefficient identities. Our contribution is an algebraically explicit experimental study that connects learned representations, decomposition growth and candidate limiting behaviour, identifying concrete objects and questions for exact follow-up.

\Cref{sec: preliminaries} recalls the required background. \Cref{sec: total tensor,sec: network of MatMul} develop the network construction, and~\Cref{sec: MatMul via neural learning} presents the finite-rank experiments and their geometric diagnostics.~\Cref{sec: border rank} introduces the Laurent extension and its experiments and~\Cref{sec: 4x4 experiments} tests the method on \(4\times4\) multiplication.~\Cref{sec: discussion} discusses the conclusions and limitations; implementation details are given in the appendix.

\section{Theoretical Background}
\label{sec: preliminaries}

\subsection{Basic notation and concepts}

We use a field \(\mathbb K\) for the exact algebraic constructions. When the ground field matters, we write \(R_{\mathbb K}(T)\) for tensor rank over \(\mathbb K\). If \(T\) has rational coefficients, scalar extension gives
\[
R_{\mathbb C}(T_{\mathbb C})
\leq
R_{\mathbb R}(T_{\mathbb R})
\leq
R_{\mathbb Q}(T).
\]

For border rank, \(\underline R_{\mathbb C}\) denotes the usual complex border rank, defined by Zariski closure, whereas \(\underline R_{\mathbb R}\) denotes real border rank, defined by Euclidean closure of the set of real tensors of rank at most \(r\). For a real tensor,
\[
\underline R_{\mathbb C}(T_{\mathbb C})
\leq
\underline R_{\mathbb R}(T).
\]

All learnable parameters in our numerical experiments are real. Consequently, the negative numerical evidence and the working hypotheses obtained from those experiments concern real rank and real border rank. On the other hand, any exact real decomposition or exact real degeneration would also give the corresponding upper bound over \(\mathbb C\).

\begin{definition}
Let \(V_1,\ldots,V_d\) be finite-dimensional vector spaces over \(\mathbb K\). The rank \(R_{\mathbb K}(T)\) of
\(T\in V_1\otimes\cdots\otimes V_d\) is the smallest integer \(r\) for which
\[
T=\sum_{s=1}^{r}
u_s^{(1)}\otimes\cdots\otimes u_s^{(d)},
\qquad
u_s^{(i)}\in V_i.
\]
When the ground field is understood, we simply write \(R(T)\).

Over \(\mathbb C\), the border rank \(\underline R_{\mathbb C}(T)\) is the smallest integer \(r\) such that \(T\) belongs to the Zariski closure of the set of tensors of rank at most \(r\). Over \(\mathbb R\), the real border rank \(\underline R_{\mathbb R}(T)\) is defined analogously using Euclidean closure in the real tensor space.
\end{definition}

The approach we present is based on the Bhattacharya--Mesner product \cite{MESNER1990204,MESNER1994595}, a $d$-ary operation in the space of $d$-order tensors that generalises matrix multiplication.

\begin{definition}
    \label{def: BMP} Let $T_{1}, T_{2}, \dots, T_{d}$ be order-$d$ tensors of types $l \times n_{2} \times \dots \times n_{d}$, $n_{1} \times l \times n_{3}\times \dots \times n_{d}$, $\cdots, n_{1} \times \dots \times l$, respectively. The \textit{generalised Bhattacharya--Mesner product} (\textit{BMP}), denoted as $\circ (T_{1}, T_{2}, \dots, T_{d})$, is a $d$-ary operation whose result is a tensor $T$ of size $n_{1}\times n_{2}\times \cdots \times n_{d}$ such that:
    \[
        (T)_{i_1, \dots, i_d} = \sum \limits_{h=1}^{l}(T_{1})_{h, i_2, \dots,i_d} \cdots (T_{d})_{i_1, i_2, \dots, h}.
    \]
\end{definition}

\subsection{Matrix Multiplication}

Let $U:=\operatorname{Mat}_{a\times b}(\mathbb K)$, $V:=\operatorname{Mat}_{b\times c}(\mathbb K)$, $W:=\operatorname{Mat}_{a\times c}(\mathbb K)$.
Matrix multiplication is the bilinear map
\[
\mu_{a,b,c}:U\times V\longrightarrow W,
\qquad
(A,B)\longmapsto AB.
\]
The associated tensor is
\[
M_{a,b,c}
=
\sum_{i=1}^{a}\sum_{j=1}^{b}\sum_{k=1}^{c}
E_{i, j}^{*}\otimes E_{j, k}^{*}\otimes E_{i, k}
\in U^{*}\otimes V^{*}\otimes W,
\]
where \(E_{i, j}^{*}\) denotes the coordinate functional dual to the matrix unit \(E_{ij}\). We also write
\[
M_{a,b,c}=\langle a,b,c\rangle.
\]

Throughout the paper we use row-major vectorisation: $\operatorname{vec}_{\mathrm r}(A)
=
(a_{11},\ldots,a_{1b},
 a_{21},\ldots,a_{ab})^{\mathsf T}$.
 
Using the standard coordinate identifications, let
\[
H=[h_1\ \cdots\ h_r]\in\mathbb K^{ab\times r},
\qquad
K=[k_1\ \cdots\ k_r]\in\mathbb K^{bc\times r},
\]
and let the \(s\)-th row of
\(F\in\mathbb K^{r\times ac}\) be \(f_s^{\mathsf T}\).
Then
\[
T=\sum_{s=1}^{r}h_s\otimes k_s\otimes f_s
\]
acts on \(A,B\) according to
\[
\operatorname{vec}_{\mathrm r}(T(A,B))
=
F^{\mathsf T}
\left[
\bigl(H^{\mathsf T}\operatorname{vec}_{\mathrm r}(A)\bigr)
\star
\bigl(K^{\mathsf T}\operatorname{vec}_{\mathrm r}(B)\bigr)
\right].
\]
Thus every choice of \(H,K,F\) defines a bilinear map of tensor rank at most \(r\). If this map agrees identically with matrix multiplication, then
\[
R(\langle a,b,c\rangle)\leq r.
\] A base-$k$ algorithm with $r$ scalar multiplications can be recursed on blocks: multiplying two $N\times N$ matrices using $k\times k$ blocks yields the recurrence
\[
    T(N)=r\,T(N/k)+O(N^{2}),
\]
because the linear pre/post-processing costs $O(N^{2})$. Now $T(N) = O \!\big(N^{\log_k r}\big)$, hence the optimal exponent $\omega$ of matrix multiplication satisfies
\[
    \omega \;=\; \inf_{k\ge 2}\,\log_{k} R(\langle k, k, k\rangle).
\]

Strassen's algorithm~\cite{strassen1969gaussian} multiplies two \(2\times2\) matrices using \(7\) bilinear multiplications instead of the usual \(8\), establishing \(R(\langle2,2,2\rangle)\leq7\). Applied recursively to blocks, it yields an algorithm with arithmetic complexity \(O(n^{\log_2 7})\). Winograd and Hopcroft--Kerr proved that six multiplications do not suffice, and hence 
\[
R(\langle2,2,2\rangle)=7
\]
\cite{winograd1971multiplication,hopcroft1971minimizing}.
De Groote subsequently classified optimal bilinear algorithms for \(2\times2\) matrix multiplication, establishing their essential equivalence to Strassen's algorithm~\cite{degroote1978optimal}.
Landsberg~\cite{landsberg2006border} later proved the corresponding equality for border rank,
\[
\underline R(\langle2,2,2\rangle)=7.
\]

\subsection{Artificial Neural Networks}

\begin{definition}
Let $V_{0}, V_{L}$ be finite-dimensional vector spaces. A \textit{feed-forward artificial neural network} is a parametrised family of functions 
\[
    f_{\theta}: V_{0} \longrightarrow V_{L},
\]
obtained by composing affine transformations and fixed nonlinear operators. Given intermediate spaces $V_{1}, \dots, V_{L-1}$, the network introduces internal variables
\[
    h_{0} = x \in V_{0}, \qquad h_{\ell} = \sigma_{\ell}\!\big( W_{\ell}(h_{\ell-1}) \big), \quad \ell = 1, \dots, L,
\]
where each $W_{\ell}: V_{\ell-1}\to V_{\ell}$ is an affine map and $\sigma_{\ell}: V_{\ell} \to V_{\ell}$ is a prescribed nonlinear activation. The elements $h_{1}, \dots, h_{L-1}$ are the \emph{hidden states} of the network and encode intermediate representations of the input. The parameter set $\theta = (W_{1}, \dots, W_{L})$ determines the realised function $f_{\theta}$.
\end{definition}

Given a distribution $\mathcal{D} \subset V_{0} \times V_{L}$ of input-output pairs $(x, y)$, learning consists in solving the optimisation problem
\[
    \min_{\theta}\; \mathbb{E}_{(x,y)\sim\mathcal{D}}\big[\,\mathcal{L}( f_{\theta}
    (x),\, y )\,\big],
\]
where $\mathcal{L}: V_{L} \times V_{L}\to \mathbb{R}_{\ge 0}$ is a prescribed loss function. In practice, the expectation is approximated by a sample average over a finite training set, and $\theta$ is updated iteratively by an optimisation algorithm at a prescribed learning rate, with one full pass through the training set called an \emph{epoch}.

In this work, we specialise this general framework to an architecture whose layers and intermediate spaces are explicitly aligned with the multilinear structure of matrix multiplication: in \textsc{MatMulNet}, the parameters play the role of tensors in a candidate rank-$r$ decomposition, and learning amounts to fitting these tensors so that the network realises the desired bilinear map.

\section{Total tensor of a network}
\label{sec: total tensor} 

In~\cite{chiantini2025producttensorsdescriptionnetworks} the total tensor of a network is computed algebraically from the topology of the underlying graph and from the activation tensors.

Let $G$ be a \emph{Directed Acyclic Graph} (DAG) with set of \emph{vertices}, or \emph{nodes}, or \emph{neurons}, $V(G)$, and set of \emph{edges} $E(G)$. A node $b \in V(G)$ is called a \emph{parent} of another node $a \in V(G)$ if there exists a directed edge $(b, a) \in E(G)$; we write
\[
    \text{parents}(a) = \{ b \in V(G) : (b, a) \in E(G) \}.
\]

The $\text{in-degree}$ of $a$, denoted by $\text{in-deg}(a)$, is $\vert \text{parents} (a) \vert$. The \emph{cardinality} of $G$ is $\vert V(G) \vert$.

\begin{definition}
    \label{def: network} In our setting, a \textit{(tensor neural) network} $\mathcal{N} = \mathcal{N}
    (G, \mathcal{A})$ is defined by the following series of data:
    \begin{enumerate}
        \item[(\textit{i})] a directed graph $G = (V(G), E(G))$ without (oriented) cycles, whose nodes represent the discrete variables of the network;  \label{enum: 1}

        \item[(\textit{ii})] a set $\mathcal{A}$ of activation tensors $A_{i}$, one for each node $a_{i}\in V(G)$, which determines the output given by the node, based on the inputs it receives;

        \item[(\textit{iii})] an ordering of the nodes which is compatible with the partial ordering given by the graph (topological sort).
            \label{enum: 3}
    \end{enumerate}
\end{definition}

Since $G$ is acyclic (\Cref{enum: 1}$(i)$), a total ordering of the nodes compatible with the partial ordering defined by the graph always exists, so the requirement in~\Cref{enum: 3}$(iii)$ can be met.

Networks defined by underlying directed graphs are useful for studying how a signal propagates from the \emph{sources} of the network, nodes with $\text{ in-deg}(a_{i}) = 0$, to its \emph{sinks}, nodes with no outgoing arrows. The variable associated with the vertex $v_{i}\in V(G)$, corresponding to the signal received and emitted by the node, assumes values in an $n_{i}$-element state set, identified with $\{0, 1, \dots, n_{i}- 1\}$. A node is said to be \emph{visible} if its state is directly observed, otherwise it is \emph{hidden}. Moreover, if the variables of all the nodes are associated to the same state set $\{0, 1, \dots, n - 1 \}$, the network is $n$-ary.

The activation tensors of~\Cref{def: network}$(ii)$ are defined as follows.

\begin{definition}
    \label{def: activation} Let $a_{i}$ be a node of $G$ with $d = \text{\it in-deg}(a_{i})$. An \emph{activation tensor} $A_{i}\in \mathcal{A}$ is any tensor of order $d + 1$ of type $n \times \cdots \times n$, such that the entry $(A_{i})_{i_1, \dots, i_d, i_{d+1}}$ represents the \emph{frequency} with which the node $a$ assumes the state $i_{d+1}$ after receiving the signal $i_{1}$ from $b_{1}$, the signal $i_{2}$ from $b_{2}$, $\dots$, and the signal $i_{d}$ from $b_{d}$.
\end{definition}

In particular, if $d = 0$, then an activation tensor at $a_{i}$ is a vector of $\bbC^{n}$. If $d = 1$, then an activation tensor at $a_{i}$ is a $n\times n$ matrix of complex numbers. In the binary case, we assume that node $a_{i}$ is inactive when its output is $0$ and $a_{i}$ is activated if its output is $1$.

\Cref{def: activation} is stated in the general probabilistic setting in which activation tensors encode transition frequencies among discrete states. In this work, however, we consider networks in which the activation tensors take on arbitrary values, without any normalisation or probabilistic interpretation.

The final state of the network is a string indexed by the nodes in the network and can be described by introducing the concept of \emph{total tensor}, whose computation is introduced in the following definition.

\begin{definition} \label{def: total tensor}
    Consider a network $\mathcal{N}(G, \mathcal{A})$ with $q$ nodes $a_1, \dots, a_q$ and states $n_1, \dots, n_q$, respectively.  For $j = 1, \dots, q$, let $A_j \in \mathcal{A}$ be the activation tensor of node $a_j$, and let $P_j = (k : a_k \in \mathrm{parents}(a_j))$ be the tuple of parent indices of $a_j$, ordered according to the total ordering of the network. Given a multi-index $(i_1, \dots, i_q)$, we write $i_{P_j} := (i_{k_1}, \dots, i_{k_{d_j}})$ for the tuple of state indices of the parents of $a_j$, where $P_j = (k_1, \dots, k_{d_j})$, and $(i_{P_j},\, i_j)$ for the ordered tuple obtained by appending the output index $i_j$ at the end. The \emph{total tensor} is the order-$q$ tensor $N$ of type $n_1 \times \cdots \times n_q$ defined by  
    $$N_{i_1, \dots, i_q} := \prod_{j = 1}^q (A_j)_{(i_{P_j},\, i_j)},$$
    so that $(A_j)_{(i_{P_j},\, i_j)}$ is the entry of $A_j$ at the multi-index whose first $d_j$ coordinates are the state indices of the parents of $a_j$ (in the order fixed by $P_j$) and whose last coordinate is the state of $a_j$ itself.
\end{definition}

The total tensor can also be computed as the BMP of the node activation tensors, appropriately resized so that each of them has order equal to the number of nodes in the graph~\cite{chiantini2025producttensorsdescriptionnetworks}).

\subsubsection{Tensor Operations}

The first two operations that resize tensors are blow and forget: blow raises the order of a tensor by one, forget raises it from $d$ to any prescribed $d' > d$.

\begin{definition} \label{def: blow}
    Let $T \in (\mathbb{K}^n)^{\otimes d}$. The \emph{blow}, or \emph{inflation}, of $T$ is the order-$(d+1)$ tensor $\mathrm{blow}(T) \in (\bbK^n)^{\otimes (d+1)}$ defined by
    \begin{equation*}
        (\text{blow}(T))_{i_1, \dots, i_{d+1}}= \bigg \{
        \begin{array}{ll}
            T_{i_1, \dots, i_d} & \text{if} \ i_1 = i_{d+1}; \\
            0                   & \text{else.}               \\
        \end{array}
    \end{equation*}
\end{definition}

\begin{definition} \label{def: forget}
    Let $T \in (\bbK^n)^{\otimes d}$ and let $d' \geq d$. For $J \subset \{1,\dots,d'\}$ with $|J| = d'- d$, the $J$-th \emph{forget}, or $J$-th \emph{copy}, of $T$ is the tensor $f_J(T) \in (\bbK^n)^{\otimes d'}$ of order $d'$ defined by
    \begin{equation*}
        (f_{J}(T))_{i_1, \dots, i_{d'}}= T_{j_1, \dots,j_{d}}
    \end{equation*}
    where $(j_{1}, \dots,j_{d})$ is the list obtained from $(i_{1}, \dots, i_{d'}
    )$ by erasing the element $i_{k}$ for all $k\in J$.
    In particular, when \(d'=d\), we have \(J=\varnothing\) and
\(f_{\varnothing}(T)=T\).
\end{definition}

To compute the final state of a network that includes hidden nodes, one performs a contraction operation on the network’s total tensor with respect to the indices of the hidden nodes.

\begin{definition}
\label{def: contraction}
Let
$T\in
\mathbb K^{a_1}\otimes\cdots\otimes\mathbb K^{a_d}$
and let \(J\subseteq[d]:=\{1,\ldots,d\}\). Write $[d]\setminus J=\{j_1<\cdots<j_q\}$.
The \emph{contraction of \(T\) along the modes in \(J\)} is the order-\(q\)
tensor \(\operatorname{contr}_{J}(T)\) defined by
\[
\bigl(\operatorname{contr}_{J}(T)\bigr)_{k_1,\ldots,k_q}
=
\sum_{\substack{1\leq i_\ell\leq a_\ell\\ \ell\in J}}
T_{i_1,\ldots,i_d},
\qquad
i_{j_s}=k_s\quad(s=1,\ldots,q).
\]
Thus contraction sums over the indices in \(J\) and retains the
indices in \([d]\setminus J\).
\end{definition}

\begin{definition} \label{def: marginal total tensor}
Let $\mathcal{N}(G,\mathcal{A})$ be a network with total tensor $N$ and let $H \subset \{1,\dots,q\}$ be the set of indices corresponding to hidden nodes. The \emph{marginal total tensor} of $\mathcal{N}$ is the tensor obtained by contracting, according to~\Cref{def: contraction}, the total tensor $N$ along the indices in $H$:
$$N^{\mathrm{obs}}_{i_{j_1},\dots,i_{j_m}} = \sum_{i_k, k \in H} N_{i_1,\dots,i_q},$$
where $j_1,\dots,j_m$ are the indices of the observed nodes.
\end{definition}

The order in which the tensor operations are performed is fixed by the topology of the graph and by the chosen total ordering, so the construction is deterministic. Consider a network $\mathcal{N}(G, \mathcal{A})$ with $\left | V(G) \right | = q$ nodes and total ordering $(a = a_{1}, a_{2},\dots, a_{q} = b)$, with \(a_1\) a source and \(a_q\) a sink of the network. For any node $a_{i}\in V(G)$, we define $P_{i}$ as the set of indexes of the parents of $a_{i}$ (see~\Cref{def: total tensor}), so that $P_{i}\subset\{1, \dots, i-1\}$.
Let \(p_i\) be the cardinality of \(P_i\) and let \(A_i\) be the activation tensor of node \(a_i\). According to
\Cref{def: activation}, \(A_i\) has order \(p_i+1\).

\begin{procedure} \label{procedure} 
    For any node $a_j \in V(G)$, we construct the associated \emph{augmented activation tensor} $T_j \in (\bbK^n)^{\otimes q}$ as follows.
    
    \begin{enumerate}
        \item [\textit{Step 1.}] If $j > 1$, apply the forget operation to $A_j$ along the indices of the predecessors of $a_j$ in the total order, that are not its parents:
        $$A'_j = f_J(A_j), \qquad J = \{1,\dots,j-1\} \setminus P_j.$$
        If $j = 1$, set $A'_j = A_j$.
        \item[\textit{Step 2.}] If $j \neq q$, apply the blow operation to $A'_j$:
        $$A''_j = \mathrm{blow}(A'_j).$$
        If $j = q$, set $A''_j = A'_j$.
        \item[\textit{Step 3.}] Apply the forget operation to $A''_j$ along the remaining indices, such that:
        $$T_j = f_H(A''_j), \qquad H = \{j+2,\dots,q\}.$$       
    \end{enumerate}
\end{procedure}

\Cref{procedure} provides a systematic construction that associates to each node $a_j \in V(G)$ an order-$q$ augmented activation tensor $T_j$, compatible with the fixed total ordering of the nodes.

Finally, the total tensor of the network $\mathcal{N}$ can be computed as
\[
    N = \circ (T_{q}, T_{1}\dots, T_{q-1}),
\]
where $T_{1}, \dots, T_{q}$ are the tensors obtained by applying~\Cref{procedure} to the activation tensors $A_{1}, \dots, A_{q}\in \mathcal{A}$ (see~\cite{chiantini2025producttensorsdescriptionnetworks}).

\section{Networks of matrix multiplication}
\label{sec: network of MatMul}

We now apply the total-tensor construction of~\Cref{sec: total tensor} to the classical and Strassen algorithms
for \(2\times2\) matrix multiplication. For each algorithm, we specify the underlying graph and the activation tensors, and verify that contracting the total tensor gives the matrix product. The Strassen construction also provides the template for the learnable architecture introduced in~\Cref{sec: MatMul via neural learning}.

The computation proceeds in three stages: we form linear combinations of the entries of each input matrix, multiply
them in pairs, one factor from each input, and recombine the resulting products linearly. With \(r\) such products,
the associated tensor is a sum of at most \(r\) rank-one terms.

\subsection{Tensor network for the $\mathbf{2 \times 2}$ matrix multiplication}
\label{subsec: 2x2 tensor network} 

Let \begin{equation} \label{eq: matrices A and B}
    A =
\begin{pmatrix}
    a_{11} & a_{12} \\
    a_{21} & a_{22}
\end{pmatrix}, \quad B =
\begin{pmatrix}
    b_{11} & b_{12} \\
    b_{21} & b_{22}
\end{pmatrix} \quad
\in \mathbb{C}^{2}\otimes \mathbb{C}^{2}.
\end{equation} 
The matrix multiplication $A \cdot B$ is a matrix $C =
\begin{pmatrix}
    c_{11} & c_{12} \\
    c_{21} & c_{22}
\end{pmatrix}
\in \mathbb{C}^{2}\otimes \mathbb{C}^{2}$ such that
\begin{equation}
    \label{eq: matrix multiplication} C =
    \begin{pmatrix}
        a_{11}b_{11}+ a_{12}b_{21} & a_{11}b_{12}+ a_{12}b_{22} \\
        a_{21}b_{11}+ a_{22}b_{21} & a_{21}b_{12}+ a_{22}b_{22}
    \end{pmatrix}.
\end{equation}

This product can be computed as the marginal total tensor of the network whose underlying graph, shown in~\Cref{fig: nonStrassen}, has three nodes $N_{1}, N_{2}, N_{3}$, in the only possible total ordering consistent with the partial order, and node $N_{2}$ is hidden. Consider a binary network and let $A_{1} = (1, 1)$ be the initial distribution of the source node $N_{1}$. The activation tensors of the nodes $N_{2}$ and $N_{3}$ are the two matrices $A$ and $B$, respectively, that is
\[
    A_{2} = A =
    \begin{pmatrix}
        a_{11} & a_{12} \\
        a_{21} & a_{22}
    \end{pmatrix}, \qquad A_{3} = B =
    \begin{pmatrix}
        b_{11} & b_{12} \\
        b_{21} & b_{22}
    \end{pmatrix}.
\]

\begin{figure}[]
    \centering
    \begin{tikzpicture}[
        >=stealth, 
        shorten >=1pt, 
        auto,
        node distance=3cm, 
        semithick 
    ]
        \tikzstyle{every state}
        =[ draw = black, thick, fill = white, minimum size = 12mm ]
        \node[state] (0) {$N_{1}$};
        \node[state] (1) [below of=0] {$N_{2}$};
        \node[state] (2) [right of=1] {$N_{3}$};
        \path[->] (0) edge node {} (1);
        \path[->] (1) edge node {} (2);
    \end{tikzpicture}
    \caption{Underlying graph of the tensor neural network whose marginal total tensor equals the matrix product. The topology of this graph ensures that the multiplication is carried out via the classical algorithm. The only possible total ordering of the nodes consistent with the partial one is $(N_1, N_2, N_3)$, and node $N_2$ is hidden.}
    \label{fig: nonStrassen}
\end{figure}

According to~\Cref{procedure}, the order of the activation tensors has to be increased based on the topology of the underlying graph of~\Cref{fig: nonStrassen}. Node $N_{1}$ transmits the signal to node $N_{2}$, which, in turn, transmits it to $N_{3}$.
Therefore, the augmented activation tensor of
$N_{1}$ becomes
\[
    T^{(1)}= f_{3}(\text{blow}(A_{1})).
\]

Node $N_{2}$ transmits the signal to $N_3$, and its augmented activation tensor is computed with a blow operation. Node $N_{3}$, on the other hand, receives the signal from its parent $N_2$, which in turn received it from $N_1$, so the activation tensor $A_{3}$ undergoes a forget operation on the first index. Therefore,
\[
    T^{(2)} = \text{blow}(A_{2}), \qquad T^{(3)}= f_{1}(A_{3}).
\]

The total tensor of the three-node network is
\[
N=\circ\bigl(T^{(3)},T^{(1)},T^{(2)}\bigr),
\qquad
N_{ijk}=a_{ij}b_{jk},
\qquad i,j,k\in\{1,2\}.
\]
Since \(N_2\) is hidden, we sum over its index to obtain the marginal total tensor, as in~\Cref{def: marginal total tensor}:
\[
N^{\mathrm{obs}}
=
\operatorname{contr}_{\{2\}}(N)
=
\begin{pmatrix}
a_{11}b_{11}+a_{12}b_{21}
&
a_{11}b_{12}+a_{12}b_{22}
\\
a_{21}b_{11}+a_{22}b_{21}
&
a_{21}b_{12}+a_{22}b_{22}
\end{pmatrix}
=AB.
\]

\subsubsection{From Strassen's algorithm to Strassen's tensor neural network} \label{subsubsec: strassen algorithm} 

Strassen's algorithm computes matrix $C$ by performing only $7$ multiplications. In fact, the entries $c_{ij}$, with $i, j = 1, 2$, depend on $7$ intermediate terms:
\begin{equation}
    \begin{aligned}
        p_{1} & = (a_{11} + a_{22}) \cdot (b_{11} + b_{22})  \\
        p_{2} & = (a_{21} + a_{22}) \cdot (b_{11})          \\
        p_{3} & = a_{11} \cdot (b_{12} - b_{22})             \\
        p_{4} & = a_{22} \cdot (b_{21} - b_{11})             \\
        p_{5} & = (a_{11} + a_{12}) \cdot b_{22}            \\
        p_{6} & = (a_{21} - a_{11}) \cdot (b_{11} + b_{12})  \\
        p_{7} & = (a_{12} - a_{22}) \cdot (b_{21} + b_{22}),
    \end{aligned}\label{eq: Strassen intermediates}
\end{equation}
such that
\[
    \begin{pmatrix}
        c_{11} & c_{12} \\
        c_{21} & c_{22}
    \end{pmatrix}
    =
    \begin{pmatrix}
        p_{1}+ p_{4}- p_{5}+ p_{7} & p_{3}+ p_{5}               \\
        p_{2}+ p_{4}               & p_{1}- p_{2}+ p_{3}+ p_{6}
    \end{pmatrix}.
\]

To compute the product of two matrices $A$ and $B$ of~\Cref{eq: matrices A and B} using Strassen’s algorithm, we consider a tensor neural network with five nodes arranged in three layers: an input layer consisting of the sources $N_A$ and $N_B$, a hidden layer consisting of $N_H$ and $N_K$, and a sink node $N_F$. The topology of the graph underlying the network is shown in~\Cref{fig: Strassen} and the initial distributions of the sources $N_{A}$ and $N_{B}$ are $D_{A} = (a_{11}, a_{12}, a_{21}, a_{22}, 0, 0, 0)$ and $D_{B }= (b_{11}, b_{12}
, b_{21}, b_{22}, 0, 0, 0)$, respectively. 

\begin{figure}[!ht]
    \centering
    \begin{tikzpicture}[
        >=stealth, 
        shorten >=1pt, 
        auto,
        node distance=3cm, 
        semithick 
    ]
        \tikzstyle{every state}
        =[ draw = black, thick, fill = white, minimum size = 12mm ]
        \node[state] (0) at (0, 2.0) {$N_{A}$}; 
        \node[state] (1) at (0, 0) {$N_{H}$}; 
        \node[state] (2) at (3, 2.0) {$N_{B}$}; 
        \node[state] (3) at (3, 0) {$N_{K}$}; 
        \node[state]
            (4)
            at
            (1.5, -1.5)
            {$N_{F}$}; 
        \path[->] (0) edge node {} (1);
        \path[->] (2) edge node {} (3); 
        \path[->] (1) edge node {} (4); 
        \path[->] (3) edge node {} (4); 
    \end{tikzpicture}
    \caption{Underlying graph of the tensor neural network whose output is obtained by contracting the total tensor
    and realises Strassen's algorithm for matrix multiplication. The network consists of five nodes $N_A, N_B, N_H, N_K, N_F$, nodes $N_H$ and $N_K$ being hidden.}
    \label{fig: Strassen}
\end{figure}

The activation tensors of the nodes $N_{H}$ and $N_{K}$ are constructed on the basis of the intermediate products $p_1, \dots, p_7$ of Strassen’s algorithm in~\Cref{eq: Strassen intermediates}. In particular, the activation tensor $H$ of node $N_H$ determines the linear combinations in matrix $A$, and the activation tensor $K$ of node $N_{K}$ determines the linear combinations in matrix $B$. The activation tensor $F$ of $N_{F}$, on the other hand, determines how the intermediate products $p_{1}, \dots, p_{7}$ are combined to obtain the product matrix $C$. Therefore, both $H$ and $K$ are $7 \times 7$ matrices:
\[
    H =
    \begin{pmatrix}
        1 & 0 & 1 & 0 & 1 & -1 & 0  \\
        0 & 0 & 0 & 0 & 1 & 0  & 1  \\
        0 & 1 & 0 & 0 & 0 & 1  & 0  \\
        1 & 1 & 0 & 1 & 0 & 0  & -1 \\
        0 & 0 & 0 & 0 & 0 & 0  & 0  \\
        0 & 0 & 0 & 0 & 0 & 0  & 0  \\
        0 & 0 & 0 & 0 & 0 & 0  & 0
    \end{pmatrix}, \qquad K =
    \begin{pmatrix}
        1 & 1 & 0  & -1 & 0 & 1 & 0 \\
        0 & 0 & 1  & 0  & 0 & 1 & 0 \\
        0 & 0 & 0  & 1  & 0 & 0 & 1 \\
        1 & 0 & -1 & 0  & 1 & 0 & 1 \\
        0 & 0 & 0  & 0  & 0 & 0 & 0 \\
        0 & 0 & 0  & 0  & 0 & 0 & 0 \\
        0 & 0 & 0  & 0  & 0 & 0 & 0
    \end{pmatrix}.
\]

The activation of \(N_{F}\), by contrast, is an order-\(3\) tensor \(\mathcal{F}\) obtained from the matrix
\[
    F_{0}=
    \begin{pmatrix}
        1  & 0  & 0 & 1 & 0 & 0 & 0 \\
        0  & 0  & 1 & -1 & 0 & 0 & 0 \\
        0  & 1  & 0 & 1 & 0 & 0 & 0 \\
        1  & 0 & 1 & 0 & 0 & 0 & 0 \\
        -1 & 1  & 0 & 0 & 0 & 0 & 0 \\
        0  & 0  & 0 & 1 & 0 & 0 & 0 \\
        1  & 0  & 0 & 0 & 0 & 0 & 0
    \end{pmatrix}.
\]
More precisely, we define \[ \mathcal{F}_{i,j,k} = \delta_{i,j}(F_{0})_{i,k}, \qquad i,j,k\in\{1,\ldots,7\}, \] where \(\delta_{i,j}\) denotes the Kronecker delta. The indices \(i\) and \(j\) correspond to the two factors entering the products \(p_{1},\ldots,p_{7}\). The condition \(i=j\) ensures that only the two factors carrying the same product index are multiplied, while the \(k\)-th column of \(F_{0}\) determines how the resulting products are combined into the \(k\)-th output coordinate.
To construct the total tensor of the network, we first restrict it to the two subgraphs of~\Cref{fig: subnetwork a b}: the one on the nodes \(N_A\), \(N_H\)~(\Cref{fig: subnetwork a}) and the one on the nodes \(N_B\), \(N_K\)~(\Cref{fig: subnetwork b}), each carrying the corresponding activation tensors. For the first restriction, the only possible total ordering is \((N_{A}, N_{H})\).
Node \(N_{A}\) transmits the signal to node \(N_{H}\), and its activation vector undergoes a blow operation:
\[
A_{1} = \text{blow}(D_{A}).
\]
By similar reasoning,
\[
B_{1} = \text{blow}(D_{B}).
\]

\begin{figure}[]
    \centering
    \begin{subfigure}[b]{0.45\textwidth}
        \centering
        \begin{tikzpicture}[
            >=stealth, 
            shorten >=1pt, 
            auto,
            node distance=3cm, 
            semithick 
        ]
            \tikzstyle{every state}
            =[ draw = black, thick, fill = white, minimum size = 12mm ]
            \node[state] (1) at (3, 2) {$N_{A}$};
            \node[state] (2) at (3, 0) {$N_{H}$};
            \path[->] (1) edge node {} (2);
        \end{tikzpicture}
        \caption{}
        \label{fig: subnetwork a}
    \end{subfigure}
    \hfill
    \begin{subfigure}[b]{0.45\textwidth}
        \centering
        \begin{tikzpicture}[
            >=stealth, 
            shorten >=1pt,
            auto,
            node distance=3cm, 
            semithick 
        ]
            \tikzstyle{every state}
            =[ draw = black, thick, fill = white, minimum size = 12mm ]
            \node[state] (1) at (1, 2) {$N_{B}$};
            \node[state] (2) at (1, 0) {$N_{K}$};
            \path[->] (1) edge node {} (2);
        \end{tikzpicture}
        \caption{}
        \label{fig: subnetwork b}
    \end{subfigure}
    \caption{Restriction of the tensor neural network with underlying graphs consisting of nodes $N_A$, $N_H$ (\ref{fig: subnetwork a}) and $N_B$, $N_K$ (\ref{fig: subnetwork b}).}
    \label{fig: subnetwork a b}
\end{figure}

According to the convention for the BMP of~\Cref{def: BMP}, the BMP of two order-\(2\) tensors is \(\circ(X,Y) = Y X\). Hence, the total tensors of the two tensor restrictions are
\[
    S_{0}^{(1)} = \circ(H,A_{1}) = A_1 H,
    \qquad
    S_{0}^{(2)} = \circ (K,B_{1}) = B_1 K.
\]

To compute the output at \(N_F\) for the prescribed source vectors, we sum over the other four node indices. We perform these contractions in two stages. 

The two modes of \(S_0^{(1)}\) correspond, in order, to \(N_A\) and \(N_H\), and those of \(S_0^{(2)}\) to
\(N_B\) and \(N_K\). We first contract the source indices, retaining the indices of \(N_H\) and \(N_K\):
\[
S^{(1)}
=
\operatorname{contr}_{\{1\}}\bigl(S_0^{(1)}\bigr),
\qquad
S^{(2)}
=
\operatorname{contr}_{\{1\}}\bigl(S_0^{(2)}\bigr).
\]
Writing \(d_A:=D_A^{\mathsf T}\) and \(d_B:=D_B^{\mathsf T}\), this gives
\[
S_s^{(1)}
=
\sum_{i=1}^{7}(d_A)_i H_{is},
\qquad
S_s^{(2)}
=
\sum_{j=1}^{7}(d_B)_j K_{js},
\qquad s=1,\ldots,7,
\]
or, equivalently,
\[
S^{(1)}=H^{\mathsf T}d_A,
\qquad
S^{(2)}=K^{\mathsf T}d_B.
\]
Thus, \(S_s^{(1)}\) and \(S_s^{(2)}\) are, respectively, the first and second factors of the product \(p_s\).

We now use \(S^{(1)}\) and \(S^{(2)}\) as the activation vectors of the two sources \(N_{S_1}\) and \(N_{S_2}\) in the three-node network of~\Cref{fig: hidden Strassen}, whose sink is \(N_F\). These two sources retain the state indices of \(N_H\) and \(N_K\). Their indices will be contracted in the second stage, after coupling them through \(\mathcal F\).

\begin{figure}[]
    \centering
    \begin{tikzpicture}[
        >=stealth, 
        shorten >=1pt, 
        auto,
        node distance=3cm, 
        semithick 
    ]
        \tikzstyle{every state}
        =[ draw = black, thick, fill = white, minimum size = 12mm ]
        \node[state] (1) at (0, 0) {$N_{S_1}$}; 
        \node[state] (3) at (3, 0) {$N_{S_2}$}; 
        \node[state]
            (4)
            at
            (1.5, -1.5)
            {$N_{F}$}; 
        \path[->] (1) edge node {} (4); 
        \path[->] (3) edge node {} (4);
    \end{tikzpicture}
    \caption{Final sub-network which computes the Strassen algorithm.}
    \label{fig: hidden Strassen}
\end{figure}

According to signal transmission, the activation tensors of node \(N_{S_1}\) and \(N_{S_2}\) become, respectively,
\[
    S_{f}^{(1)} = f_{3}\!\left(\text{blow}\!\left(S^{(1)}\right)\right), \qquad
    S_{f}^{(2)} =    \text{blow}\!\left(f_{1}\!\left(S^{(2)}\right)\right).      
\]

The source nodes \(N_{S_1}\) and \(N_{S_2}\) represent the intermediate states retained from \(N_H\) and \(N_K\).
Contracting their indices completes the output evaluation, leaving only the coordinate associated with \(N_F\).

Following~\Cref{def: marginal total tensor}, the activation tensor associated with the sink \(N_F\) appears first in the Bhattacharya--Mesner product. 
Hence, the final output is
\[
    T
    =
    \text{contraction}^{\{1,2\}}
    \left(
        \circ
        \left(
            \mathcal{F},
            S_{f}^{(1)},
            S_{f}^{(2)}
        \right)
    \right).
\]     
Indeed,~\Cref{def: BMP} gives \[ \left[ \circ \left( \mathcal{F}, S_{f}^{(1)}, S_{f}^{(2)} \right) \right]_{i,j,k} = \delta_{i,j}(F_{0})_{i,k} S_{i}^{(1)}S_{j}^{(2)}. \] Contracting the two hidden indices therefore yields \[ T_{k} = \sum_{s=1}^{7} (F_{0})_{s,k} S_{s}^{(1)}S_{s}^{(2)}, \qquad k=1,\ldots,7. \] Equivalently, \[ T = F_{0}^{\mathsf T} \left( \bigl(H^{\mathsf T}d_{A}\bigr) \star \bigl(K^{\mathsf T}d_{B}\bigr) \right), \] where \(\star\) denotes the elementwise product. Hence, \(T\) is the vector \[ T = \bigl( a_{11}b_{11}+a_{12}b_{21}, a_{11}b_{12}+a_{12}b_{22}, a_{21}b_{11}+a_{22}b_{21}, a_{21}b_{12}+a_{22}b_{22}, 0,0,0 \bigr)^{\mathsf T}. \] This expression is the padded version of the forward map used in the learning experiments. In fact, the first four columns of \(F_{0}\) form the post-combination matrix \(F\in\mathbb{R}^{7\times4}\), while the last three zero columns encode only the padding required by the present \(7\)-state tensor-network construction.

\section{Matrix Multiplication via Neural Learning}
\label{sec: MatMul via neural learning}

This section presents a neural architecture whose learnable parameters are the three matrices $H$, $K$ and $F$ characterising a matrix multiplication algorithm. The network performs all the operations described in~\Cref{subsubsec: strassen algorithm}, with matrices $H$, $K$ and $F$ initialised randomly (or, in part of the experiments of~\Cref{sec: 4x4 experiments}, as perturbations of a known scheme) and subsequently updated using a gradient descent mechanism. From now on, $n$ denotes the dimension of the matrices to be multiplied, that is $n \times n$ matrices, whilst $r$ represents the number of rank-one terms the architecture provides to decompose the matrix multiplication tensor.

\subsection{Rank experimental setup} \label{subsec: rank experimental setup}

\subsubsection{Dataset}
The training and validation sets consist of $10,000$ and $1,000$ triples respectively, composed of two random matrices $A$ and $B$ with entries sampled uniformly from $[-1, 1]$, together with their product $AB$, which corresponds to the target value. Thus, the inputs and the target can be initialised as the triple $(A, B, AB)$. 
Both sets are drawn once, before training, and kept fixed; the training set is reshuffled at each epoch.

\subsubsection{Network Architecture}
\label{subsubsec: network architecture}

The \textsc{MatMulNet} architecture has three learnable parameter matrices: $H \in \mathbb{R}^{n^2 \times r}$, $K \in \mathbb{R}^{n^2 \times r}$, and $F \in \mathbb{R}^{r \times n^2}$. Their entries are initialised independently from $\mathcal{N}(0, \alpha^{2})$, with $\alpha = 1.0$.

The inputs to the model are the flattened versions $a, b \in \mathbb{R}^{n^2}$ of $A, B \in \mathbb{R}^{n \times n}$. The target value of matrix multiplication is $v_{\text{true}} = \text{vec}(AB)$, while the predicted one $v_{\text{pred}}$ is computed by the network with the procedure shown in~\Cref{subsubsec: strassen algorithm}, using the parameter matrices $H, K$ and $F$ of the network.

The architecture is independent of the matrix size: the underlying tensor neural network (the graph in~\Cref{fig: Strassen}) and the $H/ K/ F$ parameterisation are exactly the same for $2\times 2$, $3 \times 3$ and $4 \times 4$ multiplication, and only the dimensions of the learnable parameters change with $n$. The same architecture applies, in principle, to $\langle n, n, n\rangle$ for any~$n$. We test its behaviour at larger scale on $4\times4$ multiplication in~\Cref{sec: 4x4 experiments}.

\subsubsection{Training protocol}
\label{subsubsec: rank training}

We train the network by minimising the mean squared error (MSE) between the predicted and the true flattened products over the training set $\mathbb{D}$,
\begin{equation*}
    \mathcal{L}= \frac{1}{n^2|\mathbb{D}|}\sum_{(A,B) \in \mathbb{D}}\| v_{\text{pred}}
    (A, B) - v_{\text{true}}(A, B) \|_{2}^{2},
\end{equation*}
where the factor $n^{2}$ in the denominator normalises by the number of output components, so that $\mathcal{L}$ represents the average squared error per entry of the predicted product.

For the optimisation phase we use the Adam optimiser~\cite{kingma2014adam} with mini-batch gradient descent (batch size $32$), with an initial learning rate $\eta = 10^{-4}$, halved on plateaus of the validation loss, up to a minimum value of $10^{-6}$. To prevent exploding gradients during tensor operations, we incorporate gradient clipping, which rescales the gradient whenever its norm exceeds a prescribed threshold (set to $10.0$). We report the best validation loss achieved during training for each run. 
The $10^{6}$ epochs of~\Cref*{tab: hyperparameters} are the programmed maximum: each run was stopped manually, with no quantitative stopping threshold fixed in advance, and the epochs it completed are listed in~\Cref*{tab: stop rank}.
The completed epoch counts differ across configurations, so the reported comparisons describe the searches actually performed, rather than performance at a fixed number of training epochs.

The model and the training loop are implemented in Python using the PyTorch framework. The $4\times4$ experiments replace the empirical loss with its expectation and use a compiled implementation of the same Adam update; see~\Cref{sec: 4x4 experiments}. A detailed pseudocode listing and a link to the source code are provided in~\Cref{app: code}.

\subsubsection{Calibration on \(2\times2\) multiplication} \label{subsubsec: calibration 22}

Before turning to the \(3\times3\) sweep, we apply the same finite-rank parametrisation to \(2\times2\) matrix multiplication, for which the exact rank is known to be \(7\). The matrices \(H\), \(K\), and \(F\) are again randomly initialised and learned only from triples \((A, B, AB)\); the coefficients of Strassen's decomposition are not supplied to the optimiser. This experiment is intended to check that the training procedure can reach a known finite-rank regime. Learning a \(7\)-multiplication model from examples has already been demonstrated with other multiplier-network formulations ~\cite{elser2016network,tschannen2018strassennets}. Here the experiment serves as a calibration of the specific BMP-derived architecture and Adam training protocol used in the subsequent \(3\times3\) experiments. Since Strassen's displayed factors have integer coefficients, the exact benchmark lies inside the real parameter space searched by the optimiser.

We consider $r = 6$ and $r = 7$, with $7$ independent runs each: at $r = 6$ no exact decomposition exists.

\begin{table}[h!]
    \centering
    \caption{Summary statistics of the best validation loss for rank $6$ and $7$ of $2 \times 2$ matrix multiplication ($7$ runs per configuration), together with the number of runs converging below the threshold $\tau = 10^{-12}$.}
    \label{tab: rank 6 7 results}
    \begin{tabular}{c c c c c c}
        \toprule $r$ & mean & std & median & min & successes ($< \tau$) \\
        \midrule
        $6$ & $2.716 \times 10^{-02}$ & $1.638 \times 10^{-03}$ & $2.721 \times 10^{-02}$ & $2.493 \times 10^{-02}$ & $0/7$  \\
        $7$ & $1.080 \times 10^{-19}$ & $2.190 \times 10^{-19}$ & $7.484 \times 10^{-21}$ & $1.279 \times 10^{-22}$ & $7/7$  \\
        \bottomrule
    \end{tabular}
\end{table}

\Cref{tab: rank 6 7 results} summarises the statistics of the best validation losses for all $7$ independent runs. The best validation losses of the two configurations do not overlap. At $r = 6$, every run plateaus in the narrow interval $[2.49 \times 10^{-2},\,2.96 \times 10^{-2}]$ and none of the $7$ runs reaches the threshold $\tau = 10^{-12}$. At $r = 7$, all $7$ runs fall below the threshold, with losses ranging from $1.279 \times 10^{-22}$ to $5.924 \times 10^{-19}$.

The gap between the two regimes spans approximately seventeen orders of magnitude, far exceeding the run-to-run variability observed within either configuration. This is consistent with the known equality
\[
R(\langle2,2,2\rangle)=7
\]
\cite{strassen1969gaussian,winograd1971multiplication}.
At \(r=6\), the parametrisation is algebraically incapable of representing matrix multiplication exactly, and the seven runs cluster in a narrow loss interval. At \(r=7\), by contrast, every run drives the validation MSE below the prescribed threshold, with no overlap between the two configurations.

This complete separation shows that the training protocol reaches highly accurate numerical \(7\)-term bilinear models from random initialisations. Since the quantities in~\Cref{tab: rank 6 7 results} are finite-sample validation losses, rather than full-tensor residuals, this calibration alone does not certify that the learned factors satisfy the coefficient identities of matrix multiplication exactly. Nor does it show that the learned scheme coincides with the explicit factors of Strassen's algorithm in~\Cref{subsubsec: strassen algorithm}.

\subsubsection{Experiments on $3 \times 3$ multiplication} \label{subsubsec: 3x3 experiments}

We probe the real rank of the \(3\times3\) matrix-multiplication tensor by varying the summand budget \(r\) of the model, which upper-bounds the number of scalar multiplications in the induced bilinear algorithm. The relevant characteristic-zero bounds are
\[
19
\leq
R_{\mathbb C}(M_3)
\leq
R_{\mathbb R}(M_3)
\leq
R_{\mathbb Q}(M_3)
\leq
23.
\]
The upper bound is realised by Laderman's integer-coefficient scheme~\cite{laderman1976}; it is therefore contained in the real parameter space searched by our network. Bläser's lower bound implies that \(r<19\) cannot yield an exact decomposition over either \(\mathbb R\) or \(\mathbb C\)~\cite{blaser2003}.

At the same time,
\[
17
\leq
\underline R_{\mathbb C}(M_3)
\leq
\underline R_{\mathbb R}(M_3)
\leq
20.
\]
The complex lower bound is due to Conner--Harper--Landsberg, while Smirnov's functional length-\(20\) solution is defined over any field~\cite{conner-harper-landsberg-2023,smirnov}. Thus, at \(r=17,18\), an exact finite decomposition is provably impossible, but a border-like approximation is not ruled out. At \(r=16\), even the complex border-rank lower bound excludes convergence to \(M_3\) through length-\(16\) decompositions. Moreover, because Smirnov's length-\(20\) degeneration is available over \(\mathbb R\), failure of the optimiser to access it cannot be attributed to a mismatch between the ground field of the known construction and that of the network.

We therefore sweep \(r\in\{16,\ldots,24\}\), covering the rank and border-rank windows discussed above. We include \(r=24\) to compare the optimisation behaviour with that at \(r=23\), where an exact scheme is already known. Any length-\(23\) representation can be embedded in the model with \(24\) terms by adding a zero summand. Thus every fit representable at \(r=23\) is also representable at \(r=24\), but independent training runs need not attain the same loss or success rate.
The comparison tests how the additional summand affects the numerical search. All remaining hyperparameters are given in~\Cref{tab: hyperparameters}. For each value of \(r\), we perform \(7\) independent runs.

\begin{table}[h!]
    \centering
    \caption{Summary of training hyperparameters for the rank experiments.}
    \label{tab: hyperparameters}
    \begin{tabular}{l c}
        \toprule
        \textbf{Hyperparameter} & \textbf{Value} \\
        \midrule
        Matrix size & $n = 3$ (i.e., $3 \times 3$ matrices) \\
        Rank & $r \in \{16,17,18,19,20,21,22,23,24\}$ \\
        Training samples & $10{,}000$ 
        \\
        Validation samples & $1{,}000$ \\
        Batch size & $32$ \\
        Epochs & $10^{6}$ scheduled
        \\
        Learning rate & $10^{-4}$, halved on plateau, down to $10^{-6}$
        \\
        Value range in matrices & $[-1, 1]$ \\
        \bottomrule
    \end{tabular}
\end{table}

\subsection{Results}
\label{subsec: results}

\begin{figure}[h!]
    \centering
    \begin{subfigure}[b]{0.48\textwidth}
        \centering
        \includegraphics[width=\textwidth]{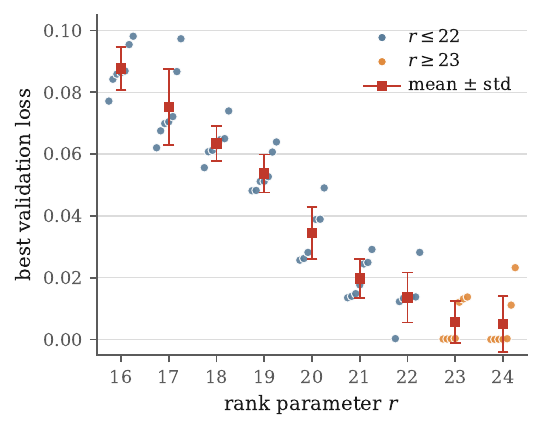}
        \caption{Linear scale}
        \label{fig: rank strip linear}
    \end{subfigure}
    \hfill
    \begin{subfigure}[b]{0.48\textwidth}
        \centering
        \includegraphics[width=\textwidth]{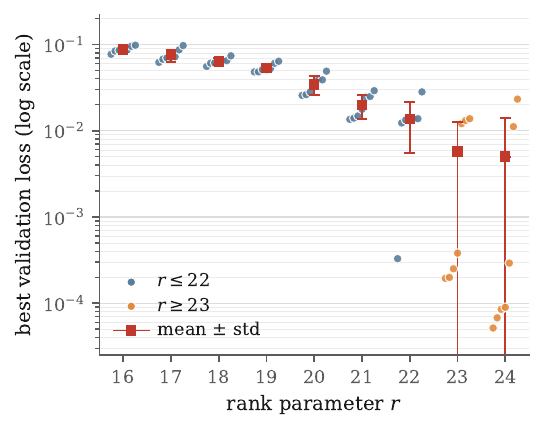}
        \caption{Logarithmic scale}
        \label{fig: rank strip log}
    \end{subfigure}
    \caption{Best validation loss for each independent run in the rank sweep $r\in\{16,\dots,24\}$, shown on a linear scale \textbf{(a)} and a logarithmic scale \textbf{(b)}.
    Runs at $r=23$ and $r=24$ are coloured in orange, since a rank-$23$ decomposition is already known to exist~\cite{laderman1976}; the red square and error bar show the mean $\pm$ one standard deviation.
    The linear scale highlights the steady decrease in mean loss from $r=16$ to $r=23$; the logarithmic scale reveals the empty band separating the converged runs from the rest, and the order-of-magnitude drop in best-run loss between $r=21$ and $r=22$.}
    \label{fig: rank strip}
\end{figure}

\Cref{tab: rank results} reports the best validation MSE and the corresponding full-tensor relative residual $\rho$, while~\Cref{fig: rank strip} shows the distribution of the validation MSE across runs.

\begin{table}[h!]
    \centering
   \caption{Summary statistics of the best validation MSE and of the full-tensor relative residual $\rho=\|T-M_{3}\|_{F}/\|M_{3}\|_{F}$ of the learned tensor, for the rank sweep ($7$ runs per configuration).} \label{tab: rank results}
    \small
    \begin{tabular}{c c c c c c c}
        \toprule
        & \multicolumn{4}{c}{validation loss} & \multicolumn{2}{c}{residual $\rho$} \\
        \cmidrule(lr){2-5}\cmidrule(lr){6-7}
        $r$ & mean & std & median & min & median & min \\
        \midrule
        $16$ & $8.77\times10^{-2}$ & $7.05\times10^{-3}$ & $8.64\times10^{-2}$ & $7.71\times10^{-2}$ & $0.513$ & $0.480$ \\
        $17$ & $7.52\times10^{-2}$ & $1.23\times10^{-2}$ & $7.05\times10^{-2}$ & $6.21\times10^{-2}$ & $0.474$ & $0.441$ \\
        $18$ & $6.34\times10^{-2}$ & $5.61\times10^{-3}$ & $6.27\times10^{-2}$ & $5.56\times10^{-2}$ & $0.434$ & $0.398$ \\
        $19$ & $5.37\times10^{-2}$ & $6.15\times10^{-3}$ & $5.12\times10^{-2}$ & $4.82\times10^{-2}$ & $0.392$ & $0.388$ \\
        $20$ & $3.45\times10^{-2}$ & $8.53\times10^{-3}$ & $3.43\times10^{-2}$ & $2.57\times10^{-2}$ & $0.335$ & $0.278$ \\
        $21$ & $1.98\times10^{-2}$ & $6.29\times10^{-3}$ & $1.77\times10^{-2}$ & $1.36\times10^{-2}$ & $0.228$ & $0.200$ \\
        $22$ & $1.36\times10^{-2}$ & $8.08\times10^{-3}$ & $1.36\times10^{-2}$ & $3.28\times10^{-4}$ & $0.196$ & $0.031$ \\
        $23$ & $5.72\times10^{-3}$ & $6.84\times10^{-3}$ & $3.80\times10^{-4}$ & $1.94\times10^{-4}$ & $0.033$ & $0.024$ \\
        $24$ & $5.00\times10^{-3}$ & $9.05\times10^{-3}$ & $8.95\times10^{-5}$ & $5.15\times10^{-5}$ & $0.017$ & $0.012$ \\
        \bottomrule
    \end{tabular}
\end{table}

For each run, the learned factors define the bilinear tensor
\[
(T_{\theta})_{ijk}
=
\sum_{s=1}^{r}H_{is}K_{js}F_{sk}.
\]
We evaluate its discrepancy from the matrix-multiplication tensor in two computational ways: on sampled matrix pairs, through the validation MSE, and directly on the tensor coefficients, through the full-tensor relative residual
\[
\rho(\theta)
=
\frac{\|T_{\theta}-M_n\|_{F}}{\|M_n\|_{F}}.
\]
These quantities are not independent error measures. For fixed parameters $\theta$, since the entries of $A$ and $B$ are independent and uniformly distributed on $[-1,1]$,
\[
\mathbb{E}_{A,B}\bigl[\mathcal{L}(\theta)\bigr]
=
\frac{1}{9n^{2}}\|T_{\theta}-M_n\|_{F}^{2}
=
\frac{n}{9}\rho(\theta)^{2}.
\]
In particular, for $n=3$,
\[
\mathbb{E}_{A,B}\bigl[\mathcal{L}(\theta)\bigr]
=
\frac{1}{3}\rho(\theta)^{2}.
\]
Thus the validation MSE is a sample-based evaluation of the same tensor discrepancy whose full-coefficient evaluation is expressed by $\rho$. The latter is computed directly from all $729$ coefficients of the learned tensor and involves no sampled input matrices. The mean and minimum validation losses highlight different aspects of the numerical search. The mean decreases monotonically from \(r=16\) to \(r=23\), from \(8.77\times10^{-2}\) to \(5.72\times10^{-3}\). At \(r=24\), it is \(5.00\times10^{-3}\), close to the value at \(r=23\). This smaller change between the last two means describes the observed training results; it does not follow from the padding argument and does not identify the tensor rank.

The minimum validation loss shows a sharper change. It decreases from \(7.71\times10^{-2}\) at \(r=16\) to \(1.36\times10^{-2}\) at \(r=21\), then drops by a factor of approximately \(41\) to \(3.28\times10^{-4}\) at \(r=22\).
This is the largest consecutive drop in the sweep and marks the first appearance of a run below the threshold \(10^{-3}\).
It identifies an empirical transition in the accuracy attained by the tested protocol. The nature of the low-loss candidates is examined further through the geometric diagnostics in~\Cref{subsec: border and rank interpretation}.

\subsection{Statistical analysis}
\label{subsec: statistical analysis rank}

We assess the trend visible in~\Cref{fig: rank strip} at three levels of granularity: a comparison of success frequencies at a fixed tolerance, a global monotonicity test across the whole sweep, and pairwise tests between consecutive ranks.

We classify a run as successful if its best validation loss falls below \(10^{-3}\). The observed classification is unchanged for any threshold strictly between \(3.8\times10^{-4}\), the largest successful loss, and \(1.1\times10^{-2}\), the smallest unsuccessful loss.~\Cref{tab: rank success} reports the per-rank counts.

For an exploratory comparison across the observed gap, we pool the six ranks \(r\leq21\), with \(0/42\) successes, against the three ranks \(r\geq22\), with \(10/21\) successes. Fisher's exact test~\cite{fisher1922interpretation} gives \(p=2.8\times10^{-6}\) for this grouping. This comparison quantifies the separation in success frequencies at the observed split; it is not an independent statistical validation of the location of the transition. In the tested runs, successful fits first occur at \(r=22\) and become more frequent at \(r=23,24\).

Among the unsuccessful runs, eight---four at \(r=22\), three at \(r=23\), and one at \(r=24\)---have full-tensor relative residuals in \([0.1943,0.1967]\), a spread of approximately \(1.2\%\). This clustering suggests recurrent optimisation behaviour at a similar error level across different ranks. 
However, comparable residual values do not establish that the runs reach the same local minimum or a common basin of attraction. 

At \(r\geq22\), the samples contain well-separated low-loss and higher-loss outcomes. Consequently, the per-rank mean reflects both the error levels attained and the relative frequencies of these outcomes, rather than only a gradual improvement in approximation accuracy.

Hence, over the whole sweep we test the null hypothesis of identical loss distributions across ranks against the ordered alternative (that the loss decreases with $r$) using the distribution-free Jonckheere-Terpstra test~\cite{jonckheere1954,terpstra1952} rather than a parametric trend test. For every pair of ranks $i<j$, let $U_{ij}$ count the pairs of runs, one from each rank, in which the $j$-th rank's loss is smaller than the $i$-th's (ties counted with weight $\nicefrac12$); the test statistic sums these counts over all $\binom{9}{2}$ rank pairs,
\[
J = \sum_{i<j} U_{ij}.
\]
Under the null hypothesis, \(J\) has known mean and variance.
We standardise it as
\[
Z=
\frac{J-\mathbb{E}[J]}
{\sqrt{\operatorname{Var}(J)}}
\]
and use its asymptotic normal approximation to compute the one-sided \(p\)-value. We obtain \(Z=9.63\) and \(p\approx3.0\times10^{-22}\), supporting the ordered alternative of decreasing losses across the sweep. The test does not require normality of the per-rank loss distributions; the reported \(p\)-value is nevertheless an asymptotic approximation, not an exact finite-sample probability.

The global trend and the observed separation in success frequencies do not establish whether every consecutive step \(r\to r+1\) produces a statistically detectable decrease in mean loss. We therefore run one-tailed Welch $t$-tests~\cite{welch1947generalization} between each pair of consecutive ranks. For two ranks with sample means $\bar x_1,\bar x_2$, standard deviations $s_1,s_2$ and $n_1 = n_2 = 7$ runs, 
\[
t = \frac{\bar x_1-\bar x_2}{\sqrt{s_1^2/n_1+s_2^2/n_2}},
\]
with degrees of freedom given by the Welch-Satterthwaite equation (hence the non-integer $df$ in~\Cref{tab: rank welch}); the test is one-tailed because the alternative hypothesis is that the loss decreases with $r$. Since we run $8$ such comparisons, the familywise error rate is controlled with the Holm-Bonferroni procedure~\cite{holm1979simple}: $p$-values are sorted $p_{(1)}\le\dots\le p_{(8)}$ and $H_{(i)}$ is rejected while $p_{(i)}\le \alpha/(8-i+1)$. Alongside significance we report Cohen's $d = (\bar x_1-\bar x_2)/s_{\text{pooled}}$, with $s_{\text{pooled}}^2 = (s_1^2+s_2^2)/2$, since every rank group has the same size ($n = 7$).~\Cref{tab: rank welch} reports the results. The three consecutive decreases from \(r=18\) to \(r=21\) remain significant after Holm correction (\(p_{\mathrm{Holm}}<0.03\), estimated \(d>1.6\)).
The steps \(16\to17\) and \(17\to18\) have estimated effect sizes \(d\approx1.2\), but do not reach significance after correction (\(p_{\mathrm{Holm}}\approx0.11\)). The steps \(21\to22\) and \(22\to23\) also do not reach corrected significance (\(p_{\mathrm{Holm}}=0.14\) and approximately \(0.11\), respectively).

These non-rejections do not establish equality of the mean losses, nor do they determine whether the lack of significance is due to limited power. In the transition region, comparisons of means should also be read alongside the success frequencies and the separation between low-loss and higher-loss outcomes. For \(23\to24\), the small estimated effect (\(d=0.09\), \(p_{\mathrm{Holm}}=0.43\)) is consistent with the similar observed means, but does not establish equivalence of optimisation performance.

Welch's test accommodates unequal variances but is not distribution-free. With only seven runs per rank and the separated outcomes observed at \(r\geq22\), its \(p\)-values should therefore be treated as exploratory parametric comparisons, alongside the nonparametric trend analysis and the descriptive summaries.

\begin{table}[h!]
    \centering
    \caption{Number of runs (out of $7$) converging to validation loss
    $< 10^{-3}$, by rank.}
    \label{tab: rank success}
    \begin{tabular}{c c c}
        \toprule $r$ & Successes & Rate \\
        \midrule
        $16$ & $0/7$ & $0\%$ \\
        $17$ & $0/7$ & $0\%$ \\
        $18$ & $0/7$ & $0\%$ \\
        $19$ & $0/7$ & $0\%$ \\
        $20$ & $0/7$ & $0\%$ \\
        $21$ & $0/7$ & $0\%$ \\
        $22$ & $1/7$ & $14\%$ \\
        $23$ & $4/7$ & $57\%$ \\
        $24$ & $5/7$ & $71\%$ \\
        \bottomrule
    \end{tabular}
\end{table}

\begin{table}[h!]
    \centering
    \caption{Pairwise Welch $t$-tests between consecutive ranks (one-tailed, Holm--Bonferroni corrected, $7$ runs per rank). Cohen's $d$ is computed from pooled standard deviations.}
    \label{tab: rank welch}
    \begin{tabular}{l c c c c c}
        \toprule Pair & $t$ & df & $p$ & $p_{\mathrm{Holm}}$ & $d$ \\
        \midrule
        $16$ vs $17$ & $2.34$ & $9.5$ & $2.1\!\times\!10^{-2}$ & $1.1\!\times\!10^{-1}$ & $1.25$ \\
        $17$ vs $18$ & $2.29$ & $8.4$ & $2.5\!\times\!10^{-2}$ & $1.1\!\times\!10^{-1}$ & $1.22$ \\
        $18$ vs $19$ & $3.08$ & $11.9$ & $4.8\!\times\!10^{-3}$ & $2.9\!\times\!10^{-2}$ & $1.65$ \\
        $19$ vs $20$ & $4.84$ & $10.9$ & $2.6\!\times\!10^{-4}$ & $2.1\!\times\!10^{-3}$ & $2.59$ \\
        $20$ vs $21$ & $3.66$ & $11.0$ & $1.9\!\times\!10^{-3}$ & $1.3\!\times\!10^{-2}$ & $1.96$ \\
        $21$ vs $22$ & $1.60$ & $11.3$ & $6.9\!\times\!10^{-2}$ & $1.4\!\times\!10^{-1}$ & $0.86$ \\
        $22$ vs $23$ & $1.98$ & $11.7$ & $3.6\!\times\!10^{-2}$ & $1.1\!\times\!10^{-1}$ & $1.06$ \\
        $23$ vs $24$ & $0.17$ & $11.2$ & $4.3\!\times\!10^{-1}$ & $4.3\!\times\!10^{-1}$ & $0.09$ \\
        \bottomrule
    \end{tabular}
\end{table}

\subsection{Rank and border rank interpretation} \label{subsec: border and rank interpretation}

The model used in this section has the ordinary finite-rank form $\sum_{s=1}^{r} h_s \otimes k_s \otimes f_s$.
For fixed finite parameters, it therefore searches for decompositions with at most \(r\) rank-one summands. This is the same type of parametrisation as the exact decompositions of Strassen and Laderman. Moreover, in the \(2\times2\) calibration described in~\Cref{subsubsec: calibration 22}, the same model is trained in a case where the exact rank is known to be \(7\), without supplying Strassen's coefficients to the optimiser. The explicitly \(\varepsilon\)-dependent model intended to investigate border rank is introduced separately in~\Cref{sec: border rank}.

The set of tensors of rank at most \(r\) is nevertheless not closed, so a very small loss does not by itself certify an exact rank-\(r\) decomposition. In particular, since
\[
\underline{R}(M_3)\leq 20,
\qquad M_3: = \langle3, 3, 3\rangle,
\]
the theoretical infimum of the unregularised loss is zero for every \(r \geq 20\): the finite-rank parametrisation can approach a border-rank configuration through increasingly large coefficients and cancellations.

The known border-rank upper bound ensures that arbitrarily accurate approximations exist for every \(r\geq20\), but does not ensure that the present optimiser reaches them. In the tested sweep, no run with \(r\leq21\) reaches validation loss below \(10^{-3}\), only one run does so at \(r=22\), and low-loss outcomes become more frequent at \(r=23,24\).

This difference motivates a closer comparison of the candidates found at \(r=22\) and \(r=23\). A small loss alone does not distinguish a bounded approximation to a finite decomposition from an approximation involving increasingly large and cancelling summands. The continuation and residual--weight analyses below examine this distinction for the candidates reached by our
method.

We further analysed the supplied length-\(22\) candidate by real, projected pseudo-arclength continuation. Starting from its parameter vector \(\theta_0\), we used the affine target path
\[
B(t)=M_3+t\bigl(T_{\theta_0}-M_3\bigr),
\qquad B(1)=T_{\theta_0},\quad B(0)=M_3.
\]
A fixed linear projection of the tensor equations and a fixed affine slice through \(\theta_0\) define the tracked curve; their construction and the predictor-corrector protocol are specified in~\Cref{app: continuation}. All reported tensor errors were evaluated independently on the full tensor, rather than only on the projected system.

The starting candidate had full-tensor relative residual \(\rho(\theta_0)\simeq1.37734\times10^{-2}\) and decomposition weight \(W(\theta_0)\simeq896.47\). On the tracked segment, the smallest recorded parameter was \(t_*\simeq0.410110\), before the curve turned towards increasing \(t\). At that point,
\[
\rho(\theta_*)\simeq6.53221\times10^{-3},
\qquad W(\theta_*)\simeq2191.79.
\]
The absolute residual of the projected equations was approximately \(1.48\times10^{-13}\), whereas the absolute full-tensor residual with respect to \(B(t_*)\) was \(1.72882\times10^{-2}\). Thus, satisfying the projected equations did not imply the same accuracy for the full tensor equations. The smallest recorded \(\rho\), attained at the next step, was \(6.53213\times10^{-3}\); minimising \(t\) and minimising the full-tensor residual are different criteria.

The run was stopped thirty accepted steps after the first detected turn in \(t\), without reaching \(t=0\). The decrease in the full-tensor residual accompanied by larger decomposition weight provides a local diagnostic, but does not establish divergence or certify a degeneration. It also does not exclude an exact bounded length-\(22\) decomposition elsewhere. The residual--weight refinement below provides a complementary analysis.

The gauge-invariant weight used in the following refinement has a direct geometric meaning.

\begin{lemma}\label{lem:bounded-weight}
Fix a ground field
\(\mathbb F\in\{\mathbb R,\mathbb C\}\), and let
\[
T_{\nu}
=
\sum_{s=1}^{r}
h_{s,\nu}\otimes k_{s,\nu}\otimes f_{s,\nu}
\longrightarrow T
\]
be a convergent sequence of tensors with at most \(r\) rank-one
summands. Define
\[
W_{\nu}
=
\sum_{s=1}^{r}
\|h_{s,\nu}\|\,\|k_{s,\nu}\|\,\|f_{s,\nu}\|.
\]
If the sequence \((W_{\nu})_{\nu}\) is bounded, then
\(
R_{\mathbb F}(T)\leq r.
\)
Consequently, if \(R_{\mathbb F}(T)>r\), every sequence of length-\(r\) decompositions converging to \(T\) must satisfy
\(
W_{\nu}\longrightarrow\infty.
\)
\end{lemma}

\begin{proof}
For each nonzero summand, set
\[
p_{s,\nu}
=
\|h_{s,\nu}\|\,\|k_{s,\nu}\|\,\|f_{s,\nu}\|.
\]
Using the termwise rescaling freedom, without changing the rank-one tensor, we may choose balanced representatives satisfying
\[
\|\widetilde h_{s,\nu}\|
=
\|\widetilde k_{s,\nu}\|
=
\|\widetilde f_{s,\nu}\|
=
p_{s,\nu}^{1/3}.
\]
A zero summand may be represented by three zero factors. If \(W_{\nu}\leq C\), then \(p_{s,\nu}\leq C\) for every \(s\), and hence all balanced factors belong to fixed compact balls. Since \(r\) is finite, after passing to a subsequence all factors converge. By continuity of the tensor product,
\[
T
=
\sum_{s=1}^{r}
\widetilde h_s\otimes\widetilde k_s\otimes\widetilde f_s,
\]
so \(R_{\mathbb F}(T)\leq r\).

For the final statement, if \(W_{\nu}\) did not tend to infinity, it would possess a bounded subsequence, and the preceding argument would imply \(R_{\mathbb F}(T)\leq r\).
\end{proof}

A complementary test uses the bounded-weight criterion above to compare the geometry of the best \(r=22\) and \(r=23\) solutions. For an increasing sequence of coefficient bounds \(c\), we impose
\[
|\theta_i|\leq c
\]
and apply projected Levenberg--Marquardt~\cite{kanzow2004levenberg}, warm-starting each refinement from the solution obtained at the preceding value of \(c\). Since the method is local, the values reported below are the residuals reached along the resulting warm-started branches; they are not claimed to be global minima over the corresponding boxes.

The coefficient bound is a numerical device for tracing the branches and is not invariant under termwise rescaling. We therefore compare the solutions using the full-tensor relative residual
\[
\rho(\theta)
=
\frac{\|T_{\theta}-M_3\|_F}{\|M_3\|_F}
\]
and the gauge-invariant decomposition weight
\[
W(\theta)
=
\sum_{s=1}^{r}
\|h_s\|\,\|k_s\|\,\|f_s\|.
\]

The two branches display sharply different profiles. Along the explored \(r=22\) branch, the residual is well described by the empirical power law 
\(
\rho\propto W^{-1.02}
\)
over nearly two decades, with no observed inflection towards a zero-residual solution at bounded weight. The branch reaches
\[
\rho=4.5\times10^{-6}
\qquad\text{at}\qquad
W=6.3\times10^{5}.
\]
This weight-growing behaviour is the finite-range signature expected from the divergence required by \Cref{lem:bounded-weight} if a sequence of length-\(22\) decompositions approaches \(M_3\) while \(R_{\mathbb R}(M_3)>22\).

By contrast, the same refinement applied to the best \(r=23\) run reaches
\[
\rho=5.3\times10^{-16}
\qquad\text{at}\qquad
W=76.
\]
This value lies below the estimated double-precision cancellation level \(2.2\times10^{-15}\) of that solution, so the reconstructed tensor is numerically indistinguishable from \(M_3\) in \texttt{float64}. At the comparable weight \(W=74\), the residual along the \(r=22\) branch is still \(2.0\times10^{-1}\).

Thus the refinement separates a weight-growing, border-like regime at \(r=22\) from a bounded, rank-like machine-precision regime at \(r=23\). Together with the different success frequencies in the original sweep and the pseudo-arclength continuation experiment, this contrast favours
\(
R_{\mathbb R}(M_3)=23
\)
as a working hypothesis and identifies \(r=22\) as the critical case for exact certification or exclusion. Since the refinement is local, it does not exclude a different bounded length-\(22\) decomposition elsewhere in parameter space.
The qualitative separation between \(r=22\) and \(r=23\) has a notable precedent in Elser's multiplier-network experiments~\cite{elser2016network}. Using a different training rule, Elser reported convergent rank-\(23\) runs in approximately \(64\%\) of \(1000\) trials, whereas with \(22\) multipliers he observed only slow convergence accompanied by growing weights. Our \(r=23\) success rate, \(4/7\), is numerically close but is based on a much smaller sample and should not be viewed as a quantitative replication.

The agreement is nevertheless qualitatively significant: two different optimisation procedures and network derivations display the same bounded-versus-weight-growing split at the \(22/23\) boundary. Whereas Elser monitored the largest individual coefficient, which depends on the choice of termwise scaling, our analysis uses the gauge-invariant decomposition weight
\[
W=\sum_s\|h_s\|\,\|k_s\|\,\|f_s\|
\]
and supplements it with constrained refinement and pseudo-arclength continuation. These diagnostics give the shared numerical phenomenon a direct geometric interpretation and strengthen the case for 
\[
R_{\mathbb R}(M_3)=23
\]
as a working hypothesis.
\begin{figure}[h!]
    \centering
    \includegraphics[width=0.78\textwidth]{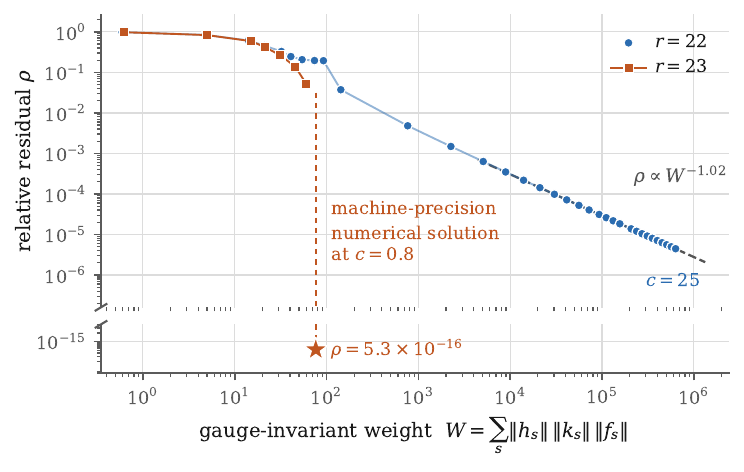}
    \caption{Local constrained refinement of the best \(r=22\) and \(r=23\) solutions. For each coefficient bound \(c\), the plotted point is the residual returned by projected Levenberg-Marquardt, warm-started from the solution obtained at the preceding value of \(c\), and is shown against the gauge-invariant decomposition weight 
    \(W=\sum_s\|h_s\|\,\|k_s\|\,\|f_s\|\).
    Along the explored \(r=22\) branch, the residual is well described by \(\rho\propto W^{-1.02}\) over nearly two decades and reaches \(4.5\times10^{-6}\) at \(W=6.3\times10^5\), without an observed transition to zero at bounded weight. At \(r=23\), the refinement reaches \(\rho=5.3\times10^{-16}\) at \(c=0.8\) and \(W=76\), below the estimated double-precision cancellation level of that solution. At the comparable weight \(W=74\), the \(r=22\) residual is still \(2.0\times10^{-1}\). The broken vertical axis separates the ordinary residual scale from the machine-precision value.}
    \label{fig: frontier}
\end{figure}

\subsubsection{On the learned decomposition}

The learned matrices $H$, $K$, and $F$ at $r = 23$ are real-valued and far from having integer entries, and the same holds after the refinement (see~\Cref{tab: decomp example} in~\Cref{app: code}). This is not unexpected. The parametrisation is non-identifiable under simultaneous permutations of the summands and termwise rescalings
\[
h_s\longmapsto \alpha_s h_s,\qquad
k_s\longmapsto \beta_s k_s,\qquad
f_s\longmapsto (\alpha_s\beta_s)^{-1}f_s,
\]
where $\alpha_s,\beta_s\neq 0$. Moreover, the loss function imposes neither sparsity nor integrality on the learned factors. The optimiser therefore has no reason to select a representative resembling Laderman's integer decomposition.

The low validation loss shows that the learned parameters encode an approximate \(23\)-term decomposition of \(\langle3,3,3\rangle\). The constrained refinement yields a bounded numerical solution of the full tensor equations with
\[
\rho
=
\frac{\|T_{\mathrm{num}}-M_3\|_F}{\|M_3\|_F}
=
5.3\times10^{-16},
\]
where every coefficient is bounded by \(0.8\) in absolute value. This residual is computed from all \(729\) coefficients of the tensor and is below the estimated cancellation level of the \texttt{float64} evaluation. Thus the refined tensor is numerically indistinguishable from \(M_3\) at double precision, although no symbolic or certified exact equality is claimed.
The numerical computation is carried out over \(\mathbb R\), not over \(\mathbb Q\). Although every stored \texttt{float64} number has a dyadic rational representation, floating-point optimisation uses rounded operations, and the resulting coefficient arrays do not by themselves constitute an exact rational decomposition. To claim a decomposition over \(\mathbb Q\), one would need to perform rational reconstruction and verify all \(729\) tensor coefficient identities exactly over \(\mathbb Q\), or after clearing denominators over \(\mathbb Z\).

A certified exact real decomposition, whether rational or not, would nevertheless also be a complex decomposition and would imply \(R_{\mathbb C}(M_3)\leq23\). Conversely, the failure of the present real optimisation to find a length-\(22\) decomposition does not exclude the existence of a complex one.

Determining whether the numerical solution is equivalent to a known exact algorithm, and possibly converting it into a more structured or human-readable representative, would require additional algebraic post-processing. Normalising the termwise scaling freedom and fixing an ordering of the summands would remove the elementary non-identifiabilities, but would not by itself produce an exact or integer-valued representative. Moreover, the large number of mutually inequivalent rank-\(23\) schemes already known~\cite{heule2021newways} makes agreement with Laderman's specific representative neither necessary nor particularly expected.

\section{Border rank}
\label{sec: border rank}

The experiments of
\Cref{subsec: rank experimental setup} use the ordinary finite-rank parametrisation and investigate the behaviour of bounded decompositions of \(M_3=\langle3,3,3\rangle\). In this section we introduce an \(\varepsilon\)-dependent extension in which all three factors are Laurent polynomials in \(\varepsilon\).

For every fixed \(\varepsilon\neq0\), the resulting tensor has rank at most \(r\). Thus the architecture searches in the natural class of one-parameter length-\(r\) families from which border-rank degenerations may arise. A genuine degeneration additionally requires the total tensor to be regular at \(\varepsilon=0\) and to have limit \(M_3\). These algebraic conditions are pursued numerically through a continuation-style schedule in which \(\varepsilon\) is progressively decreased.

Accordingly, the experiments below probe how readily gradient training produces length-\(r\) Laurent families that approximate \(M_3\) at progressively smaller positive values of \(\varepsilon\). Their purpose is to compare the accessibility of candidate degenerations at different summand budgets, rather than to provide a numerical certificate of the border rank.

The known bounds relevant to our real-valued experiments are
\[
17
\leq
\underline R_{\mathbb C}(M_3)
\leq
\underline R_{\mathbb R}(M_3)
\leq
20.
\]
The complex lower bound is due to Conner--Harper--Landsberg~\cite{conner-harper-landsberg-2023}, whereas the upper bound is realised by Smirnov's functional length-\(20\) degeneration, which is defined over the base field~\cite{smirnov,bini1980relations}. We sweep \(r\in\{16,\ldots,21\}\), thereby including one value below the known lower bound and one value above the known upper bound.

\subsection{From an $\varepsilon$-parametrisation to border rank}
\label{subsec: eps parametrisation}

A standard sufficient criterion for border rank is provided by a one-parameter degeneration.

\begin{remark}
\label{rem: eps border}
Let \(\mathbb F\in\{\mathbb R,\mathbb C\}\), let
\(T\in V_1\otimes V_2\otimes V_3\), and suppose that there exist
Laurent polynomials
\(
h_s(\varepsilon)
\in
V_1\otimes\mathbb F[\varepsilon,\varepsilon^{-1}]\), 
\(k_s(\varepsilon)
\in
V_2\otimes\mathbb F[\varepsilon,\varepsilon^{-1}]\), 
\(f_s(\varepsilon)
\in
V_3\otimes\mathbb F[\varepsilon,\varepsilon^{-1}]\), 
\(
s=1,\ldots,r\),
such that
\[
T(\varepsilon)
=
\sum_{s=1}^{r}
h_s(\varepsilon)\otimes
k_s(\varepsilon)\otimes
f_s(\varepsilon)
\]
is regular at \(\varepsilon=0\) and satisfies
\(
\lim_{\varepsilon\to0}T(\varepsilon)=T.
\)
Then
\(
\underline R_{\mathbb F}(T)\leq r.
\)
\end{remark}

For the finite Laurent expansions used in our architecture, the regularity and limit conditions can be written explicitly at the coefficient level. Let
\[
h_s(\varepsilon)=\sum_a h_{s,a}\varepsilon^a,
\qquad
k_s(\varepsilon)=\sum_b k_{s,b}\varepsilon^b,
\qquad
f_s(\varepsilon)=\sum_c f_{s,c}\varepsilon^c.
\]
Then
\[
T(\varepsilon)
=
\sum_q T_q\varepsilon^q,
\]
where
\[
T_q
=
\sum_{s=1}^{r}
\sum_{a+b+c=q}
h_{s,a}\otimes k_{s,b}\otimes f_{s,c}.
\]
The total tensor is regular at \(\varepsilon=0\) and has limit \(T\)
if and only if
\(
T_q=0
\)  for every \(q<0\),
\(
T_0=T.
\)

The architecture below parametrises the Laurent factors \(h_s(\varepsilon)\), \(k_s(\varepsilon)\), and \(f_s(\varepsilon)\). Hence, for every \(\varepsilon\neq0\), it produces a tensor of rank at most \(r\). The additional coefficient identities
\[
T_q=0\quad(q<0),
\qquad
T_0=M_3,
\]
which characterise a genuine degeneration towards \(M_3\), are not imposed identically by the present pointwise objective. Rather, the training schedule seeks them indirectly by requiring \(T(\varepsilon)\) to approximate \(M_3\) at progressively smaller positive values of \(\varepsilon\).

\subsection{$\varepsilon$-parametrised network architecture}
\label{subsec: eps architecture}

We introduce a variant of \textsc{MatMulNet} in which all three factors of the decomposition are allowed to depend on $\varepsilon$. The learnable parameters are three coefficient tensors $H_{\mathrm{coeff}}, K_{\mathrm{coeff}} \in \mathbb{R}^{n^{2} \times r \times d}$ and $F_{\mathrm{coeff}} \in \mathbb{R}^{r \times n^{2} \times d}$, where the last axis is indexed by an integer exponent $e \in [d_{\min}, d_{\max}]$ with $d = d_{\max}-d_{\min} + 1$.
All entries are initialised from $\mathcal{N}(0,\alpha^{2})$ with $\alpha = 1$.

Given \(\varepsilon>0\), the network evaluates the three Laurent factors by summing the coefficient tensors along their last axis, with weights \(\varepsilon^{d_{\min}},\ldots,\varepsilon^{d_{\max}}\), producing
\[
    H(\varepsilon)\;=\;\sum_{e=d_{\min}}^{d_{\max}}
    \varepsilon^{e}\,H_{\mathrm{coeff}}[\,\cdot\,,\,\cdot\,,e-d_{\min}],
    \qquad
    K(\varepsilon)\;=\;\sum_{e=d_{\min}}^{d_{\max}}
    \varepsilon^{e}\,K_{\mathrm{coeff}}[\,\cdot\,,\,\cdot\,,e-d_{\min}],
\]
\[
    F(\varepsilon)\;=\;\sum_{e=d_{\min}}^{d_{\max}}
    \varepsilon^{e}\,F_{\mathrm{coeff}}[\,\cdot\,,\,\cdot\,,e-d_{\min}].
\]
The matrices $H(\varepsilon)$, $K(\varepsilon)$, and $F(\varepsilon)$ are then fed into the forward pass described in~\Cref{subsubsec: rank training}, yielding the prediction
\[
    v_{\mathrm{pred}}(A, B; \varepsilon)\;=\;
    F(\varepsilon)^{\mathsf{T}}\!\Big(\big(H(\varepsilon)^{\mathsf{T}}
    \operatorname{vec}_{\mathrm r}(A)\big)\star\big(K(\varepsilon)^{\mathsf{T}}
    \operatorname{vec}_{\mathrm r}(B)\big)\Big),
\]
where \(\star\) denotes the elementwise product.

The exponent range \([d_{\min},d_{\max}]\) is a hyperparameter of the architecture. We use \([-2,4]\) as the primary configuration and compare it with the wider ranges \([-3,5]\) and \([-4,6]\) to assess sensitivity to this choice.

\subsubsection{Training protocol} \label{subsubsec: brank training protocol}
The training loss is the same per-component MSE used in~\Cref{subsubsec: rank training}, evaluated at the current value of~$\varepsilon$: 
\[
    \mathcal{L}(\varepsilon)\;=\;\frac{1}{n^2|\mathbb{D}|}
    \sum_{(A, B)\in\mathbb{D}}
    \big\|v_{\mathrm{pred}}(A, B; \varepsilon)
    -v_{\mathrm{true}}(A, B)\big\|_{2}^{2}.
\]

The parameter \(\varepsilon\) is decayed exponentially from \(\varepsilon_0=10^{-1}\) to the terminal value \(\varepsilon_{\min}=10^{-3}\), which the schedule reaches at the last of the $10^{6}$ programmed epochs. This defines a continuation-style optimisation schedule: training begins at a moderately small value of \(\varepsilon\) and progressively moves towards a regime in which the cancellations required by a degeneration become more demanding.

At each epoch, however, both \(\varepsilon\) and the learnable coefficient tensors are updated. The schedule should therefore not be identified with the evaluation of one fixed family in the limit \(\varepsilon\to0\). A small pointwise loss attained along the schedule indicates that the corresponding Laurent model approximates \(M_3\) well at that positive value of \(\varepsilon\), but does not by itself guarantee regularity of the final family at \(\varepsilon=0\).

The learning rate starts at $10^{-4}$ and is halved whenever the validation loss has not improved for $15$ consecutive epochs, until it reaches a minimum value of $10^{-6}$. The $10^{6}$ epochs are the programmed maximum, but each run was stopped manually after about one week; no quantitative stopping threshold was fixed in advance. At the interruption $\varepsilon$ had decreased only to about $0.098$-$0.099$, so the schedule was explored in a neighbourhood of its initial value $10^{-1}$ and never near its terminal value. The losses we report are the best ones observed before the interruption. The epochs completed by each run are listed in~\Cref*{tab: stop brank}.

All computations are carried out in double-precision floating-point arithmetic, to mitigate numerical errors introduced by the polynomial evaluation in $\varepsilon$, as detailed in~\Cref{subsec: numerical precision}.
The remaining hyperparameters follow the same choices described in~\Cref{subsubsec: 3x3 experiments}.
The full configuration is reported in~\Cref{tab: hyperparameters border}.

\begin{table}[h!]
    \centering
    \caption{Summary of training hyperparameters for the border-rank experiments.}
    \label{tab: hyperparameters border}
    \begin{tabular}{l c}
        \toprule
        \textbf{Hyperparameter} & \textbf{Value} \\
        \midrule
        Matrix size & $n = 3$ \\
        Rank & $r \in \{16,17,18,19,20,21\}$ \\
        Exponent range & $[-2,4]$, $[-3,5]$, $[-4,6]$ \\
        $\varepsilon$ schedule & $10^{-1} \to 10^{-3}$ (exponential decay
        over $10^{6}$ epochs
        ) \\
        Training samples & $10{,}000$ 
        \\
        Validation samples & $1{,}000$ \\
        Batch size & $32$ \\
        Epochs & $10^{6}$ programmed 
        \\
        Learning rate & $10^{-4}$\, halved on plateau, down to $10^{-6}$
        \\
        Value range in matrices & $[-1, 1]$ \\
        \bottomrule
    \end{tabular}
\end{table}

For each value of \(r\) and each exponent range, we perform \(7\) independent runs and report the minimum validation loss observed before manual termination. In these runs, \(\varepsilon\) remained close to \(10^{-1}\) throughout training. The comparison therefore concerns the initial stage of the planned continuation schedule, not its behaviour near
\(\varepsilon=0\).

\subsection{Results} \label{subsec: eps results}

\begin{figure}[h!]
    \centering
    \includegraphics[width=\textwidth]{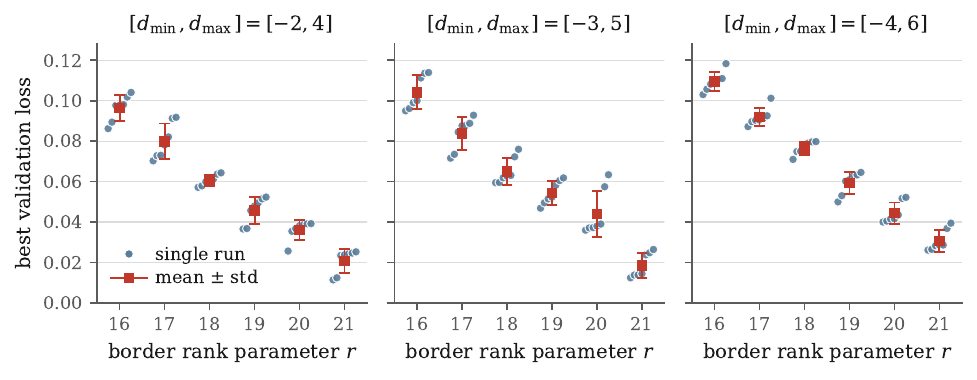}
    \caption{Distribution of best validation losses for the border-rank sweep $r\in\{16,\dots,21\}$ across three exponent ranges. Each dot is one run ($7$ per configuration); the red square and error bar show the mean $\pm$ one standard deviation.}
    \label{fig: brank strip}
\end{figure}

\Cref{tab: border results} reports the mean, standard deviation, and median of the best validation loss for the primary range $[-2, 4]$, and~\Cref{fig: brank strip} shows the distribution of individual runs across all three ranges.

\begin{table}[h!]
    \centering
    \caption{Summary statistics of the best validation loss for the border-rank sweep ($7$ runs per configuration, $[d_{\min},d_{\max}]=[-2, 4]$).}
    \label{tab: border results}
    \begin{tabular}{c c c c}
        \toprule $r$ & mean & std & median \\
        \midrule
        $16$ & $0.0964$ & $0.0064$ & $0.0978$ \\
        $17$ & $0.0799$ & $0.0088$ & $0.0781$ \\
        $18$ & $0.0607$ & $0.0027$ & $0.0607$ \\
        $19$ & $0.0457$ & $0.0066$ & $0.0479$ \\
        $20$ & $0.0361$ & $0.0048$ & $0.0382$ \\
        $21$ & $0.0207$ & $0.0061$ & $0.0237$ \\
        \bottomrule
    \end{tabular}
\end{table}

The mean validation loss decreases monotonically from $r = 16$ to $r=21$ across all three exponent ranges. Unlike the rank experiments of~\Cref{subsec: results}, no run at any tested $r$ converges below the tolerance: over the whole sweep the lowest loss is $1.1\times10^{-2}$. The distributions are unimodal and relatively tight, with no evidence of the bimodal behaviour observed in the rank case at $r\ge 22$. 

Wider exponent ranges yield systematically higher losses: the ordering $\bar x_{r,[-2,4]} < \bar x_{r,[-3,5]} < \bar x_{r,[-4,6]}$ holds at five of the six border ranks; the exception is $r=21$, where the two narrower ranges differ by less than half the run-to-run spread. This ordering is consistent with the conditioning issues discussed in~\Cref{subsec: numerical precision}, although wider exponent ranges also increase the number of parameters and may make optimisation more difficult.

\subsection{Statistical analysis} \label{subsec: statistical analysis border}

To test whether the observed decrease in validation loss with increasing $r$ is statistically significant, we apply the Jonckheere-Terpstra test~\cite{jonckheere1954,terpstra1952} for ordered alternatives. The test is one-sided, with the alternative hypothesis that losses decrease stochastically in $r$. The statistic counts the pairs of observations from different border ranks that respect the expected order; $Z$ is its excess over the null expectation, in units of the null standard deviation, and yields the $p$-value. For all three exponent ranges, the trend is highly significant (see~\Cref{tab: jt}).

\begin{table}[h!]
    \centering
    \caption{Jonckheere-Terpstra test for monotone decreasing trend in validation loss as a function of border rank $r$ ($r\in\{16,\dots,21\}$, $7$ runs per configuration). $Z$: standardised Jonckheere-Terpstra statistic; $p$: one-sided $p$-value.}
    \label{tab: jt}
    \begin{tabular}{l c c}
        \toprule $[d_{\min},d_{\max}]$ & $Z$ & $p$ \\
        \midrule
        $[-2,4]$ & $7.81$ & $< 10^{-14}$ \\
        $[-3,5]$ & $7.59$ & $< 10^{-13}$ \\
        $[-4,6]$ & $8.05$ & $< 10^{-15}$ \\
        \bottomrule
    \end{tabular}
\end{table}

We next ask whether the decreasing loss profile associated with the border-rank parameter \(r\) is robust with respect to the exponent range \([d_{\min},d_{\max}]\). We run a two-way analysis of variance~\cite{fisher1925statistical} with factors border-rank parameter and exponent range. \Cref{tab: anova} shows that both main effects are highly significant. No statistically significant interaction between the border-rank parameter and the exponent range is detected (\(p=0.29\)).
This does not establish that the effect of \(r\) is identical across exponent ranges. The decreasing profile is observed in the sample means for all three ranges and is supported by the separate trend tests in~\Cref{tab: jt}.

\begin{table}[h!]
    \centering
\caption{Two-way ANOVA (border-rank parameter \(r\) \(\times\)
exponent range) on the best validation loss
(\(r\in\{16,\dots,21\}\)).}
    \label{tab: anova}
    \begin{tabular}{l c c c c c}
        \toprule Source & SS & df & MS & $F$ & $p$ \\
        \midrule
        Rank & $0.0899$ & $5$ & $0.0180$ & $431$ & $< 10^{-15}$ \\
        Degree range & $0.0031$ & $2$ & $0.0015$ & $37.0$ & $5.8\times 10^{-13}$ \\
        Interaction & $0.0005$ & $10$ & $0.0001$ & $1.21$ & $0.29$ \\
        Residual & $0.0045$ & $108$ & $0.0000$ & & \\
        \bottomrule
    \end{tabular}
\end{table}

To test whether each consecutive drop is individually significant, we run one-tailed Welch $t$-tests between all consecutive rank pairs, with Holm-Bonferroni correction for multiple comparisons (see~\cite{welch1947generalization,holm1979simple}), on the primary exponent range $[-2,4]$. The results are shown in~\Cref{tab: brank welch}.

\begin{table}[h!]
    \centering
   \caption{Pairwise Welch \(t\)-tests between consecutive values of the
border-rank parameter \(r\) (one-tailed, Holm--Bonferroni corrected,
\(7\) runs per configuration,
\([d_{\min},d_{\max}]=[-2,4]\)).} \label{tab: brank welch}\begin{tabular}{l c c c c c}
        \toprule Pair & $t$ & df & $p$ & $p_{\mathrm{Holm}}$ & $d$ \\
        \midrule
        $16$ vs $17$ & $4.00$ & $11.0$ & $1.1\!\times\!10^{-3}$ & $2.1\!\times\!10^{-3}$ & $2.14$ \\
        $17$ vs $18$ & $5.50$ & $7.1$ & $4.3\!\times\!10^{-4}$ & $1.3\!\times\!10^{-3}$ & $2.94$ \\
        $18$ vs $19$ & $5.55$ & $7.9$ & $2.8\!\times\!10^{-4}$ & $1.1\!\times\!10^{-3}$ & $2.97$ \\
        $19$ vs $20$ & $3.11$ & $11.0$ & $4.9\!\times\!10^{-3}$ & $4.9\!\times\!10^{-3}$ & $1.66$ \\
        $20$ vs $21$ & $5.24$ & $11.4$ & $1.2\!\times\!10^{-4}$ & $6.1\!\times\!10^{-4}$ & $2.80$ \\
        \bottomrule
    \end{tabular}
\end{table}

In the primary configuration \([-2,4]\), all five consecutive decreases are significant after Holm correction, with large estimated effect sizes (\(d>1.6\)). Among the theoretically admissible consecutive boundaries \(17/18\), \(18/19\), and \(19/20\), the largest standardised separation occurs at \(18/19\), although only marginally over \(17/18\).
Within the tested architecture and training protocol, this gives a weak numerical reason to prioritise the \(18/19\) boundary for exact or certified follow-up. It does not identify the actual border rank. 

\subsection{Numerical precision analysis}
\label{subsec: numerical precision}

The \(\varepsilon\)-parametrised architecture evaluates Laurent polynomials whose monomials span a wide dynamic range. For \(0<\varepsilon<1\), the ratio between the largest and smallest monomials in each factor is
\[
\frac{\varepsilon^{d_{\min}}}{\varepsilon^{d_{\max}}}
=
\varepsilon^{d_{\min}-d_{\max}}.
\]
\Cref{tab: precision} reports this range at the scheduled terminal value \(\varepsilon=10^{-3}\).

\begin{table}[h!]
    \centering
    \caption{Dynamic range of the monomials in each Laurent factor at \(\varepsilon=10^{-3}\), before multiplication by the learned coefficients.}
    \label{tab: precision}
    \begin{tabular}{l c c c}
        \toprule
        \([d_{\min},d_{\max}]\)
        & Largest monomial
        & Smallest monomial
        & \(\log_{10}\)(ratio) \\
        \midrule
        \([-2,4]\) & \(10^{6}\)  & \(10^{-12}\) & \(18\) \\
        \([-3,5]\) & \(10^{9}\)  & \(10^{-15}\) & \(24\) \\
        \([-4,6]\) & \(10^{12}\) & \(10^{-18}\) & \(30\) \\
        \bottomrule
    \end{tabular}
\end{table}

All the displayed monomials lie well within the normal range of \texttt{float64}; their evaluation does not, by itself, cause underflow. The issue is instead relative precision. Double-precision arithmetic provides approximately \(16\) significant decimal digits. When contributions of very different magnitudes are accumulated, a small contribution may not change the rounded sum. Conversely, cancellation between large, already rounded quantities may expose and amplify earlier rounding errors~\cite{goldberg1991floating}.

The magnitudes of the actual contributions depend on the learned coefficients, not only on the powers of \(\varepsilon\).
Their numerical evaluation also depends on the order of summation and on cancellations. Consequently, the exponent range alone does not determine a fixed number of ``usable'' or ``lost'' terms, nor does it imply an automatic reduction of the polynomial degree.

The total tensor is trilinear in \(H(\varepsilon)\), \(K(\varepsilon)\), and \(F(\varepsilon)\), so its formal expansion has exponents between \(3d_{\min}\) and \(3d_{\max}\).
In a Laurent representation of a border-rank degeneration, the coefficients of negative powers cancel in the total tensor even when individual summands grow. Evaluating such a representation can therefore be sensitive to rounding errors in the factors, their products, and the subsequent summation.

The experimental results show an effect of the exponent range: wider ranges generally attain higher mean losses, and the two-way ANOVA reports a significant main effect (\(F=37.0\), \(p=5.8\times10^{-13}\);~\Cref{tab: anova}). This is consistent with poorer conditioning, but the analysis does not distinguish numerical rounding effects from optimisation difficulties associated with the larger parametrisation.

The reported minima were attained during the executed initial part of the schedule, with \(\varepsilon\) remaining close to \(10^{-1}\). At these values, the monomial ratios are approximately \(10^6\), \(10^8\), and \(10^{10}\) for the three exponent ranges. 
These ratios alone do not establish numerical stability, which also depends on the learned coefficients and cancellations.
The values in~\Cref{tab: precision} instead describe the planned endpoint \(\varepsilon=10^{-3}\), not reached by these runs, and do not by themselves explain the observed loss levels.

\section{Scaling to $4\times4$ multiplication}
\label{sec: 4x4 experiments}

We use \(M_4=\langle4,4,4\rangle\) to characterise the method at larger scale and assess the accessibility of known decompositions under controlled perturbations and random initialisation. The relevant bounds are
\[
34 \leq R(M_4) \leq 48,
\qquad
29 \leq \underline R(M_4) \leq 46.
\]
The lower bound on the rank is due to Bl\"aser~\cite{blaser2003}. The upper bound was obtained by AlphaEvolve~\cite{novikov2025alphaevolve} over \(\mathbb C\), and then with rational coefficients by Dumas, Pernet and Sedoglavic~\cite{dumas2025rational48}. For the border rank, the lower bound is due to Landsberg and Micha{\l}ek~\cite{landsberg2018lowerbound} and the upper bound to Smirnov~\cite{smirnov2014length46}. A rational length-\(48\) scheme is therefore available~\cite{dumas2025rational48}, and Strassen's algorithm applied recursively gives an exact length-\(49\) scheme, so decompositions of both lengths lie in the real parameter space searched by the network.

\subsection{Experimental setup}

\subsubsection{Protocol}
We replace the empirical loss of~\Cref{subsubsec: rank training} with its expectation over \(A,B\), derived in~\Cref{subsec: results}, computed directly on the learned tensor: 
$$\mathcal L=\|T_\theta-M_4\|_F^2/144=\tfrac49\rho^2.$$
There are no training or validation sets, and every Adam step is full-batch, with gradient clipping at \(10\). The learning rate is halved on plateaus (patience \(10^4\) steps) down to \(10^{-8}\). A run stops when \(\rho<10^{-8}\), when the best loss gains less than a relative \(10^{-3}\) over \(10^5\) steps at the minimum learning rate, or at the step budget, and ends with at most \(200\) Levenberg-Marquardt (LM) iterations. A run is classified as numerically solved if at the end \(\rho<10^{-10}\) and \(W<10^{4}\); the recovered factors need not represent the scheme used to build the start.

\subsubsection{Recovery from a perturbed known scheme}
We perturb the rational length-\(48\) scheme with Gaussian noise of standard deviation \(\sigma\) times the root mean square of its non-zero coefficients, with \(100\) starts per level. Each start is trained with initial learning rate \(10^{-2}\) and \(3\times10^{-3}\) (at most \(10^6\) steps), and is also refined by LM alone (\(1000\) iterations) as a control.

\begin{table}[h!]
    \centering
   \caption{Recovery from perturbations of a length-\(48\) scheme:
starts solved out of \(100\).}
    \label{tab: 4x4 recovery}
    \begin{tabular}{c c c c c c c}
        \toprule
        \(\sigma\) & \(0.4\) & \(0.5\) & \(0.6\) & \(0.7\) & \(0.8\) & \(1.0\) \\
        \midrule
        LM alone & \(100\) & \(100\) & \(68\) & \(0\) & \(0\) & \(0\) \\
        Training, \(\eta=10^{-2}\) & \(100\) & \(98\) & \(27\) & \(0\) & \(0\) & \(0\) \\
        Training, \(\eta=3\times10^{-3}\) & \(100\) & \(97\) & \(23\) & \(0\) & \(0\) & \(0\) \\
        \bottomrule
    \end{tabular}
\end{table}

\Cref{tab: 4x4 recovery} shows a sharp transition between \(\sigma=0.5\) and \(\sigma=0.7\). Every solved training run had already reached \(\rho<10^{-8}\) before the final LM iterations. LM alone attains the highest aggregate recovery counts, while nine starts at \(\sigma=0.6\) are solved only by training, showing that the procedures recover different subsets of the perturbed starts. This experiment quantifies the sensitivity of recovery to initialisation and the complementary outcomes of the two procedures. Eleven of the unsolved runs ended with \(\rho<0.11\). These eleven runs remained unsolved after LM refinement: \(\rho\) kept decreasing while \(W\) grew beyond \(10^4\), in five cases with \(\rho\propto W^{-1}\), the same scaling observed along the length-\(22\) branch for \(M_3\).

\subsubsection{Random starts}
We train from \(\mathcal N(0,1)\) factors at \(r\in\{46,47,48,49\}\), with initial learning rates \(3\times10^{-4}\), \(10^{-3}\), \(3\times10^{-3}\), \(10^{-2}\) and \(100\) runs each, at most \(10^7\) steps. Every run stopped on a stall.

\begin{table}[h!]
    \centering
    \caption{Random starts at \(n=4\): \(400\) runs per length, four learning rates pooled. \(\rho\) is measured after the final LM iterations, \(W\) at the end of training. \(k=64\rho^2\) at the end of training is averaged over the runs within \(5\%\) of an integer (\(1574\) of \(1600\)).}
    \label{tab: 4x4 random}
    \begin{tabular}{c c c c c c}
        \toprule
        \(r\) & solved & min \(\rho\) & median \(\rho\) & median \(W\) & mean \(k\) \\
        \midrule
        \(46\) & \(0/400\) & \(0.217\) & \(0.331\) & \(4928\) & \(7.41\) \\
        \(47\) & \(0/400\) & \(0.217\) & \(0.306\) & \(4746\) & \(6.58\) \\
        \(48\) & \(0/400\) & \(0.177\) & \(0.283\) & \(5129\) & \(5.73\) \\
        \(49\) & \(0/400\) & \(0.125\) & \(0.280\) & \(5370\) & \(5.11\) \\
        \bottomrule
    \end{tabular}
\end{table}

Nearly all runs end with \(\rho\approx\sqrt{k/64}\) for an integer \(k\geq1\). On a sample of runs, the residual has \(k\) singular values close to \(1\) in each flattening, revealing a recurrent multilinear structure in the approximation error. Direct constrained refinement of selected candidates, with bounds starting above their largest coefficient, reaches the same kind of level while \(W\) grows to \(10^5\)--\(10^6\). The sampled protocol of~\Cref{subsubsec: rank training} gave the same outcome: none of \(35\) runs at \(r=45,\dots,49\) reached a validation loss below \(10^{-3}\).

The mean \(k\) decreases regularly with \(r\), by \(0.83\), \(0.85\) and \(0.62\), with no discontinuity at \(r=48\) or \(r=49\), although a decomposition exists at both lengths. Failure from random starts is therefore not evidence of non-existence. The contrast used for \(M_3\) is informative because there the method reaches a highly accurate numerical solution at length \(r=23\), where exact decompositions are known; at \(n=4\) this positive control is not available from random starts.

Related observations of integer-valued residual plateaus and difficulties in recovering known \(4\times4\) schemes by direct optimisation were reported by Smirnov~\cite{smirnov}. Here these phenomena are examined through residual--weight diagnostics and controlled perturbation experiments.

\section{Discussion and conclusions}
\label{sec: discussion} 

In \textsc{MatMulNet}, the number \(r\) of rank-one summands is fixed by the architecture. This gives an explicit bilinear parametrisation whose learned factors can be examined as tensor decompositions. Our experiments compare both the approximation accuracy attained at different values of \(r\) and the size of the summands used to attain it.

The most informative numerical indication comes from the success frequencies together with the residual-weight profiles.
In the finite-rank sweep, no run with \(r\leq21\) attains validation loss below \(10^{-3}\), whereas this threshold is reached in \(1/7\), \(4/7\), and \(5/7\) runs at \(r=22,23,24\), respectively. These observations direct attention to the geometry of the low-loss candidates at lengths \(22\) and \(23\).

The constrained refinement sharpens this distinction. At \(r=22\), the residual decreases only together with a sustained increase of the gauge-invariant decomposition weight, following approximately \(\rho\propto W^{-1.02}\) along the explored branch. At \(r=23\), by contrast, the full-tensor equations are satisfied to machine precision at bounded weight.~\Cref{lem:bounded-weight} provides the geometric basis for this diagnostic: if a sequence of tensors represented by at most \(r\) rank-one summands converges with uniformly bounded decomposition weight, its limit has rank at most \(r\). Consequently, convergence to a tensor of rank greater than \(r\) requires the weight to diverge. Weight divergence along a particular sequence, however, does not by itself exclude a bounded length-\(r\) decomposition of the same limit.

The agreement between the sweep, the constrained refinement, and the pseudo-arclength continuation therefore favours
\[
R_{\mathbb R}(\langle3,3,3\rangle)=23
\]
as a working hypothesis and identifies \(r=22\) as the critical case for exact certification or exclusion. This is a numerical indication, not a new lower bound: the methods are local to the solutions reached and do not exclude a different exact length-\(22\) decomposition elsewhere in parameter space.

The border-rank sweep behaves differently. The mean best validation loss decreases monotonically from \(r=16\) to \(r=21\) for all three exponent ranges, with significant ordered trends in each range (\Cref{tab: jt}). In the primary configuration \([d_{\min},d_{\max}]=[-2,4]\), all five consecutive decreases remain significant after Holm correction
(\Cref{tab: brank welch}). Nevertheless, the lowest validation loss over the entire sweep is \(1.1\times10^{-2}\), and no order-of-magnitude transition comparable to that in the finite-rank experiments is observed. The ordering is consistent with \(17\leq\underline R(\langle3,3,3\rangle)\leq20\) and does not determine the border rank. 
Nevertheless, within the tested numerical protocol, the loss profile weakly prioritises \(r=19\) and identifies the \(18/19\) boundary as a natural target for subsequent exact or certified searches.
These observations concern runs that explored only a neighbourhood of the initial value \(\varepsilon=10^{-1}\); they do not describe the behaviour of the planned continuation at small \(\varepsilon\).

The \(4\times4\) experiments apply the same diagnostics to a harder case. Small residuals are not all of the same kind: solved runs attain the prescribed full-tensor residual tolerance with bounded weight, whereas other runs with a small residual keep decreasing only as \(W\) grows, in five cases as \(\rho\propto W^{-1}\), the scaling observed along the length-\(22\) branch for \(M_3\). Most unsuccessful runs instead end at \(\rho\approx\sqrt{k/64}\) for an integer \(k\geq1\), revealing recurrent residual levels among the stalled runs.

The limitation lies in the search from random starts.
It reaches no exact decomposition at lengths \(46\)--\(49\), including length \(49\), where one is known to exist. Failure of the search at a given length is therefore not evidence of non-existence. The contrast for \(M_3\) is different because the method does reach the known length \(r=23\), which gives a reference for interpreting the behaviour at \(r=22\).

\subsection{Limitations}
\label{subsec: limitations}

The weaker separation observed in the border-rank experiments may reflect both numerical conditioning and optimisation difficulty. At fixed \(\varepsilon>0\), perturbing the coefficient of \(\varepsilon^e\) in \(H\) by \(\Delta H_e\) changes the evaluated factor by
\[
\Delta H(\varepsilon)=\varepsilon^e\Delta H_e.
\]
Thus, as \(\varepsilon\) decreases, positive and negative exponents attenuate and amplify these perturbations, respectively.
The effect on the prediction additionally depends on the other two factors and on cancellations. The optimisation problem can therefore involve directions with very different sensitivities, rather than a uniformly flattening loss landscape.

Wider exponent ranges increase both this disparity of scales and the number of learnable parameters. Together with the finite-precision effects discussed in~\Cref{subsec: numerical precision}, these features may hinder the search for border-rank degenerations. The present experiments, however, do not isolate their respective contributions to the observed loss levels.

The \(4\times4\) experiments use the expected loss rather than the sampled protocol of the \(3\times3\) sweep, and the runs were stopped by a stalling criterion rather than manually, so the two cases should be compared qualitatively.

\subsection{Future work}
\label{subsec: future work}

The architecture applies unchanged to \(\langle n,n,n\rangle\) for any \(n\). Determining from which length the method reaches exact solutions from random starts at \(n=4\) remains open. The elementary non-identifiabilities of the factorisation are termwise rescalings and simultaneous permutations of the summands; they do not form a general \(\mathrm{GL}_r\)-action. Fixing these freedoms is a first step towards rational reconstruction or other canonical forms of the learned schemes.

A natural control experiment is to repeat the critical \(r=22,23\) searches with complex-valued factors. Agreement between the real and complex profiles would extend the numerical indication beyond the real parameter space, although it would still not constitute a lower bound.

Finally, the different search paradigms should be compared on common targets and summand budgets. Elser's conservative learning explores a continuous real factor space; the present method uses standard gradient training together with geometric diagnostics; StrassenNets and Stapleton introduce ternary restrictions or ternarisation; AlphaTensor searches a finite coefficient space using reinforcement learning. These methods target the same underlying tensor-decomposition problem, but not the same parameter space or notion of numerical success.

The experiments described in this work show that, when the computation graph mirrors the algebraic structure of bilinear algorithms, gradient-based training learns highly accurate seven-term models for \(2\times2\) multiplication and, combined with constrained refinement, yields a machine-precision length-\(23\) numerical solution for \(3\times3\) multiplication. The accompanying geometric diagnostics quantify the residual--weight trade-off along the explored length-\(22\) branch. At \(n=4\), they characterise recovery under perturbations and distinguish recurrent residual plateaus from weight-growing approximations, providing concrete directions for further
algebraic and numerical investigation.

\section*{Acknowledgments}

AB has been partially supported by GNSAGA of INDAM; the European Union's HORIZON-MSCA-2023-DN-JD programme under the Horizon Europe (HORIZON) Marie Skłodowska-Curie Actions, grant agreement 101120296 (TENORS); the Simons Institute; the TensorDec Laboratory of the Department of Mathematics of the University of Trento. 

SM warmly thanks Luca Chiantini and Cristiano Bocci, supervisors of her Ph.D. thesis~\cite{marziali2026}, whose results are extended in the present paper, for many helpful discussions.

{ \bibliographystyle{alpha} \bibliography{strassen} }

@article{strassen1969gaussian,
  title={Gaussian elimination is not optimal},
  author={Strassen, V.},
  journal={Numerische mathematik},
  volume={13},
  number={4},
  pages={354--356},
  year={1969},
  publisher={Springer}
}

@article{winograd1971multiplication,
  title={On multiplication of 2$\times$ 2 matrices},
  author={Winograd, S.},
  journal={Linear algebra and its applications},
  year={1971}
}

@article{hopcroft1971minimizing,
  title={On minimizing the number of multiplications necessary for matrix multiplication},
  author={Hopcroft, J.E. and Kerr, L.R.},
  journal={SIAM Journal on Applied Mathematics},
  volume={20},
  number={1},
  pages={30--36},
  year={1971},
  publisher={SIAM}
}

@article{laderman1976,
  author  = {Laderman, J.D.},
  title   = {A noncommutative algorithm for multiplying $3 \times 3$ matrices using 23 multiplications},
  journal = {Bulletin of the American Mathematical Society},
  volume  = {82},
  number  = {1},
  pages   = {126--128},
  year    = {1976},
  url     = {https://projecteuclid.org/journals/bulletin-of-the-american-mathematical-society-new-series/volume-82/issue-1/A-noncommutative-algorithm-for-multiplying-3-times-3-matrices-using/bams/1183537626.full}
}

@article{blaser2003,
  author  = {Bl{\"a}ser, M.},
  title   = {On the complexity of the multiplication of matrices of small formats},
  journal = {Journal of Complexity},
  volume  = {19},
  number  = {1},
  pages   = {43--60},
  year    = {2003},
  doi     = {10.1016/S0885-064X(02)00007-9},
  url     = {https://www.sciencedirect.com/science/article/pii/S0885064X02000079}
}

@article{conner-harper-landsberg-2023,
  author  = {Conner, A. and Harper, A. and Landsberg, J.M.},
  title   = {New lower bounds for matrix multiplication and the $3\times3$ determinant},
  journal = {Forum of Mathematics, Pi},
  volume  = {11},
  pages   = {e17},
  year    = {2023},
  doi     = {10.1017/fmp.2023.14},
  url     = {https://www.cambridge.org/core/journals/forum-of-mathematics-pi/article/new-lower-bounds-for-matrix-multiplication-and-det3/170334CE1BE6A9604BD4D941258326F3}
}

@article{fawzi2022alphatensor,
  author  = {Fawzi, A. and Balog, M. and Huang, A. and Hubert, T. and Romera-Paredes, B. and
             Barekatain, M. and Novikov, A. and Ruiz, F.J.R. and Schrittwieser, J. and
             Swirszcz, G. and Silver, D. and Hassabis, D. and Kohli, P.},
  title   = {Discovering faster matrix multiplication algorithms
             with reinforcement learning},
  journal = {Nature},
  volume  = {610},
  number  = {7930},
  pages   = {47--53},
  year    = {2022},
  doi     = {10.1038/s41586-022-05172-4}
}

@article{bini1980relations,
  title={Relations between exact and approximate bilinear algorithms. Applications},
  author={Bini, D.},
  journal={Calcolo},
  volume={17},
  number={1},
  pages={87--97},
  year={1980},
  publisher={Springer}
}

@article{schoenhage1981,
  author  = {Sch{\"o}nhage, A.},
  title   = {Partial and Total Matrix Multiplication},
  journal = {SIAM Journal on Computing},
  year    = {1981},
  volume  = {10},
  number  = {3},
  pages   = {434--455},
  doi     = {10.1137/0210031}
}

@book{bcs1997,
  author    = {B{\"u}rgisser, P. and Clausen, M. and Shokrollahi, M.A.},
  title     = {Algebraic Complexity Theory},
  publisher = {Springer},
  series    = {Grundlehren der mathematischen Wissenschaften},
  volume    = {315},
  year      = {1997},
  doi       = {10.1007/978-3-662-03338-8}
}

@article{landsberg2006border,
  title={The border rank of the multiplication of 2$\times$ 2 matrices is seven},
  author={Landsberg, J.M.},
  journal={Journal of the American Mathematical Society},
  volume={19},
  number={2},
  pages={447--459},
  year={2006}
}

@article{welch1947generalization,
  title={The generalization of ‘STUDENT'S’problem when several different population varlances are involved},
  author={Welch, B.L.},
  journal={Biometrika},
  volume={34},
  number={1-2},
  pages={28--35},
  year={1947},
  publisher={Oxford University Press}
}

@misc{chiantini2025producttensorsdescriptionnetworks,
      title={Product of Tensors and Description of Networks}, 
      author={Chiantini, L. and D'Inverno, G.A. and Marziali, S.},
      year={2025},
      eprint={2402.06768},
      archivePrefix={arXiv},
      primaryClass={math.AG},
      url={https://arxiv.org/abs/2402.06768}, 
}

@article{smirnov,
  title={The bilinear complexity and practical algorithms for matrix multiplication},
  author={Smirnov, A.V.},
  journal={Computational Mathematics and Mathematical Physics},
  volume={53},
  number={12},
  pages={1781--1795},
  year={2013},
  publisher={Springer}
}

@article{BERNARDI201778,
title = {Tensor decomposition and homotopy continuation},
journal = {Differential Geometry and its Applications},
volume = {55},
pages = {78-105},
year = {2017},
note = {Geometry and complexity theory},
issn = {0926-2245},
doi = {https://doi.org/10.1016/j.difgeo.2017.07.009},
url = {https://www.sciencedirect.com/science/article/pii/S0926224517301055},
author = {Bernardi, A. and Daleo, N.S. and Hauenstein, J.D. and Mourrain, B.}
}

@article{MESNER1990204,
  author  = {Mesner, D.M. and Bhattacharya, P.},
  title   = {Association schemes on triples and a ternary algebra},
  journal = {Journal of Combinatorial Theory, Series A},
  year    = {1990},
  volume  = {55},
  number  = {2},
  pages   = {204--234},
  doi     = {10.1016/0097-3165(90)90068-8},
  url     = {https://www.sciencedirect.com/science/article/pii/0097316590900688}
}

@article{MESNER1994595,
  author  = {Mesner, D.M. and Bhattacharya, P.},
  title   = {A ternary algebra arising from association schemes on triples},
  journal = {Journal of Algebra},
  year    = {1994},
  volume  = {164},
  number  = {3},
  pages   = {595--613},
  doi     = {10.1006/jabr.1994.1081},
  url     = {https://www.sciencedirect.com/science/article/pii/S0021869384710817}
}

@article{kingma2014adam,
  title={Adam: A method for stochastic optimization},
  author={Kingma, D.P. and Ba, J.},
  journal={arXiv preprint arXiv:1412.6980},
  year={2014}
}

@article{jonckheere1954,
  author  = {Jonckheere, A.R.},
  title   = {A distribution-free {$k$}-sample test
             against ordered alternatives},
  journal = {Biometrika},
  volume  = {41},
  number  = {1--2},
  pages   = {133--145},
  year    = {1954},
  doi     = {10.1093/biomet/41.1-2.133}
}

@book{fisher1925statistical,
  author    = {Fisher, R.A.},
  title     = {Statistical Methods for Research Workers},
  publisher = {Oliver and Boyd},
  address   = {Edinburgh},
  year      = {1925}
}

@article{fisher1922interpretation,
  title={On the interpretation of $\chi^2$ from contingency tables, and the calculation of {P}},
  author={Fisher, R.A.},
  journal={Journal of the Royal Statistical Society},
  volume={85},
  number={1},
  pages={87--94},
  year={1922},
  publisher={JSTOR}
}

@article{kanzow2004levenberg,
  author  = {Kanzow, C. and Yamashita, N. and Fukushima, M.},
  title   = {Levenberg--Marquardt methods with strong local convergence
             properties for solving nonlinear equations with convex constraints},
  journal = {Journal of Computational and Applied Mathematics},
  volume  = {172},
  number  = {2},
  pages   = {375--397},
  year    = {2004}
}

@article{holm1979simple,
  author  = {Holm, S.},
  title   = {A simple sequentially rejective multiple test procedure},
  journal = {Scandinavian Journal of Statistics},
  volume  = {6},
  number  = {2},
  pages   = {65--70},
  year    = {1979}
}

@article{terpstra1952,
  author  = {Terpstra, T.J.},
  title   = {The asymptotic normality and consistency of Kendall's test against
             trend, when ties are present in one ranking},
  journal = {Indagationes Mathematicae},
  volume  = {14},
  pages   = {327--333},
  year    = {1952}
}

@article{elser2016network,
  author  = {Elser, V.},
  title   = {A Network That Learns Strassen Multiplication},
  journal = {Journal of Machine Learning Research},
  volume  = {17},
  number  = {116},
  pages   = {1--13},
  year    = {2016}
}

@inproceedings{tschannen2018strassennets,
  author    = {Tschannen, M. and Khanna, A. and Anandkumar, A.},
  title     = {StrassenNets: Deep Learning with a Multiplication Budget},
  booktitle = {Proceedings of the 35th International Conference on
               Machine Learning},
  series    = {Proceedings of Machine Learning Research},
  volume    = {80},
  pages     = {4985--4994},
  publisher = {PMLR},
  year      = {2018}
}

@article{heule2021newways,
  author  = {Heule, M.J.H. and Kauers, M. and Seidl, M.},
  title   = {New Ways to Multiply \(3\times3\)-Matrices},
  journal = {Journal of Symbolic Computation},
  volume  = {104},
  pages   = {899--916},
  year    = {2021},
  doi     = {10.1016/j.jsc.2020.10.003}
}

@misc{stapleton2025rank23,
  author        = {Stapleton, J.},
  title         = {A 60-Addition, Rank-23 Scheme for Exact
                   \(3\times3\) Matrix Multiplication},
  year          = {2025},
  eprint        = {2508.03857},
  archivePrefix = {arXiv},
  primaryClass  = {cs.DS}
}

@article{goldberg1991floating,
  author  = {Goldberg, D.},
  title   = {What Every Computer Scientist Should Know About
             Floating-Point Arithmetic},
  journal = {ACM Computing Surveys},
  volume  = {23},
  number  = {1},
  pages   = {5--48},
  year    = {1991},
  doi     = {10.1145/103162.103163}
}

@article{degroote1978optimal,
  author  = {de Groote, H.F.},
  title   = {On varieties of optimal algorithms for the computation
             of bilinear mappings. {II}. Optimal algorithms for
             {$2 \times 2$}-matrix multiplication},
  journal = {Theoretical Computer Science},
  volume  = {7},
  number  = {2},
  pages   = {127--148},
  year    = {1978},
  doi     = {10.1016/0304-3975(78)90045-2}
}

@article{tichavsky2017numerical,
  author  = {Tichavsk{\'y}, P. and Phan, A.H. and Cichocki, A.},
  title   = {Numerical {CP} decomposition of some difficult tensors},
  journal = {Journal of Computational and Applied Mathematics},
  volume  = {317},
  pages   = {362--370},
  year    = {2017},
  doi     = {10.1016/j.cam.2016.12.007}
}

@article{daiLim2023stability,
  author  = {Dai, Z. and Lim, L.-H.},
  title   = {Numerical stability and tensor nuclear norm},
  journal = {Numerische Mathematik},
  volume  = {155},
  pages   = {345--376},
  year    = {2023},
  doi     = {10.1007/s00211-023-01377-5}
}

@article{deSilvaLim2008illposed,
  author  = {de Silva, V. and Lim, L.-H.},
  title   = {Tensor Rank and the Ill-Posedness of the Best Low-Rank
             Approximation Problem},
  journal = {SIAM Journal on Matrix Analysis and Applications},
  volume  = {30},
  number  = {3},
  pages   = {1084--1127},
  year    = {2008},
  doi     = {10.1137/06066518X}
}

@article{vermeylen2025stability,
  author  = {Vermeylen, C. and {Van Barel}, M.},
  title   = {Stability improvements for fast matrix multiplication},
  journal = {Numerical Algorithms},
  volume  = {100},
  pages   = {645--683},
  year    = {2025},
  doi     = {10.1007/s11075-024-01806-y}
}

@inproceedings{kauersMoosbauer2023flip,
  author    = {Kauers, M. and Moosbauer, J.},
  title     = {Flip Graphs for Matrix Multiplication},
  booktitle = {Proceedings of the 2023 International Symposium on
               Symbolic and Algebraic Computation},
  series    = {ISSAC '23},
  pages     = {381--388},
  publisher = {Association for Computing Machinery},
  year      = {2023},
  doi       = {10.1145/3597066.3597120}
}

@article{novikov2025alphaevolve,
  title   = {{AlphaEvolve}: A coding agent for scientific and algorithmic discovery},
  author  = {Novikov, Alexander and others},
  journal = {arXiv preprint arXiv:2506.13131},
  year    = {2025}
}

@article{dumas2025rational48,
  title   = {A non-commutative algorithm for multiplying 4x4 matrices using 48 non-complex multiplications},
  author  = {Dumas, Jean-Guillaume and Pernet, Cl{\'e}ment and Sedoglavic, Alexandre},
  journal = {arXiv preprint arXiv:2506.13242},
  year    = {2025}
}

@article{landsberg2018lowerbound,
  title   = {A $2n^2-\log_2(n)-1$ lower bound for the border rank of matrix multiplication},
  author  = {Landsberg, J. M. and Micha{\l}ek, Mateusz},
  journal = {International Mathematics Research Notices},
  volume  = {2018},
  number  = {15},
  pages   = {4722--4733},
  year    = {2018},
  doi     = {10.1093/imrn/rnx025}
}

@article{smirnov2014length46,
  title   = {The Approximate Bilinear Algorithm of Length 46 for Multiplication of 4 x 4 Matrices},
  author  = {Smirnov, A. V.},
  journal = {arXiv preprint arXiv:1412.1687},
  year    = {2014}
}

@phdthesis{marziali2026,
  title={Product of Tensors and the Algebra of Networks},
  author={Marziali, Sara},
  year={2026},
  school={University of Siena}
}

\appendix
\crefalias{section}{appendix}
\section{Training algorithm and code}
\label{app: code}

\Cref{alg: matmulnet} provides a pseudocode listing of the \textsc{MatMulNet} training loop described in~\Cref{sec: MatMul via neural learning}. The source code for the neural-training experiments, together with the code and numerical data for the projected pseudo-arclength continuation described in~\Cref{app: continuation}, is available in the same repository at \url{https://github.com/MarzialiS/MatMulNet}.

In the pseudocode, the forward pass is written in the equivalent matrix form derived from the BMP construction in~\Cref{subsubsec: strassen algorithm}. The vectors \(H^{\mathsf T}a\) and \(K^{\mathsf T}b\) have length \(r\); their elementwise product contains the \(r\) intermediate bilinear products, and multiplication by \(F^{\mathsf T}\) produces the \(n^2\) output coordinates. The loop is shown up to the programmed maximum \(E\); in the reported experiments it was interrupted manually as described
in~\Cref{subsubsec: rank training}.

\begin{algorithm}[h!]
    \caption{\textsc{MatMulNet} Training}
    \label{alg: matmulnet}
    \begin{algorithmic}[1]
        \State \textbf{Input:} rank $r$, matrix size $n$, maximum epochs $E$, batch size $B$, training-set size $N$, learning rates $\eta_{0}\to\eta_{1}$ plateau factor $\gamma$ and patience $P$
        \State \textbf{Initialize:} $H, K \sim \mathcal{N}(0, 1)^{n^2 \times r}$,
        $F \sim \mathcal{N}(0, 1)^{r \times n^2}$

        \State Initialize optimiser: Adam($\{H, K, F\}$, lr $= \eta_{0}$)

        \State Sample the training set $\mathcal{D}$ of $N$ triples and the validation set $\mathcal{V}$ once; $\ \mathcal{L}^{\star}\leftarrow \infty$

        \For{epoch $= 1$ to $E$}

        \State Reshuffle $\mathcal{D}$, whose triples are $(A, B, \operatorname{vec}_{\mathrm r}(AB))$, entries uniform in $[-1,1]$

        \For{batch $(A, B, V_{\text{true}})$ in $\mathcal{D}$}
        \State $V_{\text{pred}}\leftarrow [\ ]$

        \For{each sample $(A, B, v_t)$ in batch}
        \State $a \leftarrow \operatorname{vec}_{\mathrm r}(A)$;
               $b \leftarrow \operatorname{vec}_{\mathrm r}(B)$
        \State $S^{(1)} \leftarrow H^{\mathsf T}a$;
               $S^{(2)} \leftarrow K^{\mathsf T}b$
        \State $p \leftarrow S^{(1)} \star S^{(2)}$
        \State $v_{\text{pred}} \leftarrow F^{\mathsf T}p$
        \State Append $v_{\text{pred}}$ to $V_{\text{pred}}$
        \EndFor

        \State $\mathcal{L}\leftarrow \text{MSE}(V_{\text{pred}}, V_{\text{true}})$

        \State Backpropagate $\nabla_{H,K,F}\mathcal{L}$ \State Clip gradients:
        $\|\nabla\| \leq 10$

        \State Update parameters with Adam \EndFor

        \State $\mathcal{L}_{\text{val}}\leftarrow$ evaluate on $\mathcal{V}$

        \If{$\mathcal{L}_{\text{val}} < \mathcal{L}^{\star}$}
        \State $\mathcal{L}^{\star}\leftarrow \mathcal{L}_{\text{val}}$;
        store $(H, K, F)$ \EndIf

        \State If $\mathcal{L}^{\star}$ has not improved for $P$ consecutive epochs, multiply the learning rate by $\gamma$, without going below $\eta_{1}$

        \EndFor

        \State \textbf{Return:} the stored $H, K, F$ and the best validation
        loss $\mathcal{L}^{\star}$
    \end{algorithmic}
\end{algorithm}

\Cref{tab: decomp example} reports the numerical factors \(H\), \(K\), and \(F\) obtained from the best run at \(r=23\) after the constrained refinement of~\Cref{subsec: border and rank interpretation}, at clamp level \(c=0.8\). At full precision these factors define the tensor
\[
T_{\mathrm{num}}
=
\sum_{t=1}^{23}h_t\otimes k_t\otimes f_t
\]
and satisfy
\[
\frac{\|T_{\mathrm{num}}-M_3\|_F}{\|M_3\|_F}
=
5.3\times10^{-16}.
\]
The largest entrywise error of \(T_{\mathrm{num}}-M_3\) is \(6.7\times10^{-16}\), while the estimated \texttt{float64} cancellation level for this solution is \(2.2\times10^{-15}\).

Each summand \(t=1,\ldots,23\) is specified by the corresponding columns \(h_t,k_t\in\mathbb R^9\) of \(H,K\) and by the corresponding row \(f_t\in\mathbb R^9\) of \(F\). The entries displayed in the table are rounded to two decimal places for readability and therefore do not reproduce the reported residual. The full-precision coefficients
used for the numerical verification are provided as machine-readable supplementary data.
\begin{table}[h!]
    \centering
   \caption{Numerical \(23\)-term factors obtained by constrained refinement of the best run at \(r=23\). At full precision the resulting tensor has relative residual \(\rho=5.3\times10^{-16}\), with all coefficients bounded by \(0.8\) in absolute value. Each group of three rows lists one summand \(t\); the vectors \(h_t\), \(k_t\), and \(f_t\) are reshaped as \(3\times3\) matrices. The displayed entries are rounded to two decimal places and are not intended to reproduce the reported full-precision residual.}
    \label{tab: decomp example}
    \small
    \begin{tabular}{c|rrr|rrr|rrr}
        \hline
        $t$ & \multicolumn{3}{c|}{$h_t$} & \multicolumn{3}{c|}{$k_t$} & \multicolumn{3}{c}{$f_t$} \\
        \hline
        $1$ & $-0.44$ & $0.65$ & $0.69$ & $0.71$ & $-0.16$ & $-0.50$ & $-0.56$ & $-0.65$ & $0.41$ \\
         & $-0.38$ & $0.56$ & $0.59$ & $0.01$ & $0.00$ & $-0.01$ & $-0.57$ & $-0.66$ & $0.42$ \\
         & $-0.41$ & $0.61$ & $0.65$ & $-0.59$ & $0.13$ & $0.41$ & $-0.08$ & $-0.09$ & $0.06$ \\
        \hline
        $2$ & $0.63$ & $-0.64$ & $0.76$ & $0.40$ & $-0.36$ & $0.39$ & $0.74$ & $-0.71$ & $0.14$ \\
         & $0.17$ & $-0.17$ & $0.20$ & $-0.51$ & $0.46$ & $-0.50$ & $-0.61$ & $0.59$ & $-0.12$ \\
         & $0.44$ & $-0.45$ & $0.53$ & $0.74$ & $-0.67$ & $0.72$ & $-0.22$ & $0.22$ & $-0.04$ \\
        \hline
        $3$ & $-0.68$ & $0.08$ & $0.47$ & $0.07$ & $0.70$ & $-0.66$ & $0.21$ & $-0.52$ & $0.47$ \\
         & $0.72$ & $-0.11$ & $-0.42$ & $0.10$ & $0.72$ & $-0.70$ & $-0.25$ & $0.63$ & $-0.56$ \\
         & $-0.11$ & $0.00$ & $0.10$ & $-0.05$ & $-0.24$ & $0.25$ & $-0.21$ & $0.51$ & $-0.46$ \\
        \hline
        $4$ & $0.04$ & $-0.55$ & $-0.40$ & $-0.45$ & $0.75$ & $0.57$ & $0.72$ & $0.69$ & $0.80$ \\
         & $-0.62$ & $-0.74$ & $-0.18$ & $-0.68$ & $0.34$ & $-0.38$ & $0.80$ & $0.03$ & $0.11$ \\
         & $0.80$ & $0.36$ & $-0.18$ & $-0.20$ & $-0.34$ & $-0.80$ & $0.04$ & $0.73$ & $0.77$ \\
        \hline
        $5$ & $0.35$ & $-0.09$ & $-0.14$ & $0.75$ & $-0.18$ & $0.72$ & $-0.17$ & $0.43$ & $-0.38$ \\
         & $-0.34$ & $-0.01$ & $0.34$ & $-0.70$ & $0.15$ & $-0.71$ & $-0.12$ & $0.29$ & $-0.26$ \\
         & $0.78$ & $-0.07$ & $-0.56$ & $-0.64$ & $0.29$ & $-0.41$ & $0.29$ & $-0.73$ & $0.65$ \\
        \hline
        $6$ & $-0.06$ & $-0.80$ & $0.40$ & $-0.80$ & $-0.29$ & $0.03$ & $0.71$ & $-0.75$ & $-0.72$ \\
         & $0.70$ & $0.03$ & $0.44$ & $0.03$ & $0.60$ & $0.66$ & $-0.80$ & $0.41$ & $0.35$ \\
         & $0.76$ & $0.02$ & $0.48$ & $-0.54$ & $-0.76$ & $-0.61$ & $-0.60$ & $-0.14$ & $-0.22$ \\
        \hline
        $7$ & $0.46$ & $-0.69$ & $-0.72$ & $-0.59$ & $-0.73$ & $-0.55$ & $-0.39$ & $-0.45$ & $0.29$ \\
         & $-0.39$ & $0.59$ & $0.62$ & $0.61$ & $0.75$ & $0.57$ & $-0.26$ & $-0.31$ & $0.19$ \\
         & $-0.42$ & $0.63$ & $0.66$ & $0.08$ & $0.10$ & $0.07$ & $0.65$ & $0.76$ & $-0.48$ \\
        \hline
        $8$ & $0.72$ & $-0.10$ & $-0.46$ & $-0.74$ & $0.68$ & $0.07$ & $-0.29$ & $0.73$ & $-0.65$ \\
         & $0.19$ & $-0.03$ & $-0.12$ & $-0.01$ & $0.01$ & $0.00$ & $-0.30$ & $0.75$ & $-0.66$ \\
         & $0.51$ & $-0.07$ & $-0.32$ & $0.61$ & $-0.56$ & $-0.06$ & $-0.04$ & $0.10$ & $-0.09$ \\
        \hline
        $9$ & $0.51$ & $0.72$ & $0.62$ & $0.78$ & $-0.18$ & $-0.55$ & $0.06$ & $-0.59$ & $-0.62$ \\
         & $0.44$ & $0.62$ & $0.54$ & $0.64$ & $-0.14$ & $-0.45$ & $0.00$ & $0.02$ & $0.02$ \\
         & $0.48$ & $0.68$ & $0.59$ & $-0.11$ & $0.02$ & $0.08$ & $0.07$ & $-0.66$ & $-0.69$ \\
        \hline
        $10$ & $0.67$ & $-0.67$ & $0.80$ & $0.43$ & $0.38$ & $-0.32$ & $0.76$ & $-0.02$ & $-0.80$ \\
         & $0.17$ & $-0.18$ & $0.21$ & $-0.55$ & $-0.49$ & $0.42$ & $-0.62$ & $0.01$ & $0.66$ \\
         & $0.47$ & $-0.47$ & $0.56$ & $0.80$ & $0.71$ & $-0.60$ & $-0.23$ & $0.00$ & $0.24$ \\
        \hline
        $11$ & $-0.80$ & $-0.78$ & $-0.11$ & $-0.80$ & $-0.10$ & $0.80$ & $0.30$ & $0.80$ & $0.17$ \\
         & $0.19$ & $-0.80$ & $-0.66$ & $-0.07$ & $-0.62$ & $0.60$ & $0.35$ & $0.39$ & $-0.27$ \\
         & $-0.40$ & $-0.79$ & $-0.33$ & $0.61$ & $-0.44$ & $-0.16$ & $0.00$ & $0.50$ & $0.43$ \\
        \hline
        $12$ & $0.65$ & $-0.66$ & $0.78$ & $0.11$ & $-0.02$ & $-0.08$ & $-0.65$ & $0.39$ & $0.19$ \\
         & $0.56$ & $-0.57$ & $0.68$ & $-0.77$ & $0.17$ & $0.54$ & $0.78$ & $-0.47$ & $-0.22$ \\
         & $0.61$ & $-0.62$ & $0.73$ & $0.80$ & $-0.18$ & $-0.56$ & $0.63$ & $-0.38$ & $-0.18$ \\
        \hline
        $13$ & $-0.80$ & $-0.70$ & $-0.57$ & $0.33$ & $0.41$ & $0.31$ & $0.08$ & $-0.77$ & $-0.80$ \\
         & $0.68$ & $0.60$ & $0.48$ & $-0.43$ & $-0.52$ & $-0.40$ & $-0.07$ & $0.63$ & $0.66$ \\
         & $0.73$ & $0.64$ & $0.52$ & $0.62$ & $0.76$ & $0.57$ & $-0.02$ & $0.23$ & $0.24$ \\
        \hline
        $14$ & $0.30$ & $0.28$ & $0.04$ & $0.35$ & $-0.17$ & $0.20$ & $0.62$ & $0.60$ & $0.70$ \\
         & $-0.40$ & $-0.38$ & $-0.05$ & $0.80$ & $-0.40$ & $0.46$ & $-0.02$ & $-0.02$ & $-0.02$ \\
         & $0.79$ & $0.75$ & $0.10$ & $0.38$ & $-0.19$ & $0.22$ & $0.69$ & $0.67$ & $0.78$ \\
        \hline
        $15$ & $0.42$ & $-0.42$ & $0.50$ & $-0.15$ & $0.58$ & $0.72$ & $-0.44$ & $0.12$ & $0.32$ \\
         & $0.64$ & $-0.65$ & $0.77$ & $0.15$ & $-0.60$ & $-0.74$ & $-0.30$ & $0.08$ & $0.21$ \\
         & $-0.36$ & $0.37$ & $-0.43$ & $0.02$ & $-0.08$ & $-0.10$ & $0.74$ & $-0.21$ & $-0.53$ \\
        \hline
    \end{tabular}
\end{table}

\begin{table}[h!]
    \centering
    \caption*{~\Cref{tab: decomp example} (continued).}
    \small
    \begin{tabular}{c|rrr|rrr|rrr}
        \hline
        $t$ & \multicolumn{3}{c|}{$h_t$} & \multicolumn{3}{c|}{$k_t$} & \multicolumn{3}{c}{$f_t$} \\
        \hline
        $16$ & $-0.48$ & $0.71$ & $0.75$ & $0.03$ & $0.34$ & $-0.33$ & $0.00$ & $0.00$ & $0.00$ \\
         & $0.48$ & $-0.72$ & $-0.76$ & $0.08$ & $0.78$ & $-0.74$ & $-0.64$ & $-0.75$ & $0.47$ \\
         & $-0.09$ & $0.13$ & $0.14$ & $0.04$ & $0.37$ & $-0.35$ & $0.59$ & $0.69$ & $-0.44$ \\
        \hline
        $17$ & $0.47$ & $-0.58$ & $0.79$ & $0.05$ & $0.34$ & $-0.33$ & $0.24$ & $-0.59$ & $0.53$ \\
         & $0.64$ & $-0.54$ & $0.54$ & $-0.16$ & $-0.46$ & $0.54$ & $-0.28$ & $0.71$ & $-0.63$ \\
         & $0.54$ & $-0.56$ & $0.68$ & $0.19$ & $0.65$ & $-0.73$ & $-0.23$ & $0.58$ & $-0.52$ \\
        \hline
        $18$ & $-0.51$ & $0.37$ & $-0.31$ & $-0.29$ & $-0.19$ & $-0.69$ & $0.18$ & $-0.45$ & $0.40$ \\
         & $-0.31$ & $0.51$ & $-0.78$ & $0.36$ & $0.16$ & $0.75$ & $0.12$ & $-0.31$ & $0.27$ \\
         & $-0.14$ & $-0.24$ & $0.63$ & $-0.46$ & $0.27$ & $-0.20$ & $-0.31$ & $0.77$ & $-0.68$ \\
        \hline
        $19$ & $-0.06$ & $0.21$ & $0.80$ & $-0.21$ & $-0.75$ & $-0.75$ & $-0.51$ & $-0.16$ & $0.75$ \\
         & $0.66$ & $-0.80$ & $0.05$ & $0.31$ & $0.80$ & $0.76$ & $-0.28$ & $-0.69$ & $-0.10$ \\
         & $0.07$ & $-0.20$ & $-0.72$ & $-0.71$ & $-0.09$ & $0.21$ & $0.80$ & $0.80$ & $-0.70$ \\
        \hline
        $20$ & $0.72$ & $-0.59$ & $-0.80$ & $0.52$ & $0.52$ & $-0.45$ & $0.00$ & $-0.65$ & $0.22$ \\
         & $-0.80$ & $0.53$ & $0.80$ & $-0.32$ & $0.80$ & $-0.80$ & $-0.80$ & $-0.37$ & $0.14$ \\
         & $0.63$ & $0.36$ & $-0.09$ & $-0.70$ & $0.24$ & $-0.30$ & $0.34$ & $0.65$ & $-0.80$ \\
        \hline
        $21$ & $-0.52$ & $-0.50$ & $-0.06$ & $0.16$ & $-0.62$ & $-0.77$ & $0.66$ & $0.64$ & $0.74$ \\
         & $-0.80$ & $-0.76$ & $-0.10$ & $0.00$ & $-0.01$ & $-0.01$ & $0.67$ & $0.65$ & $0.75$ \\
         & $0.45$ & $0.43$ & $0.06$ & $-0.13$ & $0.51$ & $0.63$ & $0.09$ & $0.09$ & $0.10$ \\
        \hline
        $22$ & $0.71$ & $0.67$ & $0.09$ & $0.06$ & $0.08$ & $0.06$ & $-0.50$ & $-0.48$ & $-0.56$ \\
         & $-0.60$ & $-0.57$ & $-0.07$ & $-0.44$ & $-0.54$ & $-0.41$ & $0.60$ & $0.58$ & $0.68$ \\
         & $-0.65$ & $-0.62$ & $-0.08$ & $0.46$ & $0.56$ & $0.42$ & $0.49$ & $0.47$ & $0.55$ \\
        \hline
        $23$ & $-0.40$ & $-0.32$ & $-0.49$ & $0.10$ & $-0.40$ & $-0.49$ & $0.00$ & $0.00$ & $0.00$ \\
         & $-0.61$ & $-0.49$ & $-0.74$ & $-0.09$ & $0.34$ & $0.42$ & $0.08$ & $-0.74$ & $-0.77$ \\
         & $0.34$ & $0.28$ & $0.42$ & $-0.16$ & $0.62$ & $0.76$ & $-0.07$ & $0.69$ & $0.71$ \\
        \hline
    \end{tabular}
\end{table}

% Insert as a subsection of the computational appendix.
\subsection{Projected pseudo-arclength continuation}
\label{app: continuation}

We identify the coefficients of \(T_\theta\) with
\(f(\theta)\in\mathbb R^{729}\), using row-major ordering in each
matrix factor and storing the parameters in \(\theta\in\mathbb R^{594}\)
in factor-major, then summand-major order. The starting vector
\(\theta_0\) is read from the supplied labelled coefficient file;
no preliminary refinement is applied in this experiment. Let
\(m\in\mathbb R^{729}\) represent \(M_3\), and set
\[
b(t)=m+t\bigl(f(\theta_0)-m\bigr).
\]
The Jacobian is computed analytically. For example,
\[
\frac{\partial (T_\theta)_{abc}}{\partial (h_s)_i}
=\delta_{ai}(k_s)_b(f_s)_c,
\]
with analogous expressions for the other two factors.

Let \(J_0=L\Sigma V^{\mathsf T}\) be a singular value decomposition
of \(Df(\theta_0)\). We define its numerical rank \(q\) by the
threshold \(\sigma_i>10^{-10}\sigma_1\); for this input,
\(q=550\). Let \(Q=V_{[:,1:q]}\) and
\(S=V_{[:,q+1:594]}^{\mathsf T}\). We restrict to the fixed affine
slice by writing
\[
\theta(y)=\theta_0+Qy,\qquad y\in\mathbb R^{550},
\]
so that \(S(\theta-\theta_0)=0\). The singular vectors defining
this slice are not updated along the run.

Using NumPy's default random generator with seed \(12345\), we
sample a \(550\times729\) matrix \(R\) with independent standard
normal entries. An economy QR factorisation of \(R^{\mathsf T}\)
provides a matrix with orthonormal columns; its transpose is our
fixed projection \(P\in\mathbb R^{550\times729}\). We track the
projected system
\[
g(y,t)=\frac{P\bigl(f(\theta(y))-b(t)\bigr)}{\|M_3\|_F}=0.
\]
This projection does not make the projected and full tensor
equations equivalent. The full discrepancies are monitored separately.

Writing \(z=(y,t)\), we initialise at \(z_0=(0,1)\) and choose
the unit tangent \(\tau_0\in\ker Dg(z_0)\) with negative
\(t\)-component. Subsequent tangent vectors are obtained by a
bordered linear solve, normalised in the Euclidean norm and oriented
consistently with the previous tangent. From an accepted point
\(z_k\), the predictor is
\(z_{\mathrm{pred}}=z_k+\Delta s\,\tau_k\). A damped Newton
iteration then solves
\[
\begin{cases}
g(z)=0,\\
\tau_k^{\mathsf T}(z-z_{\mathrm{pred}})=0.
\end{cases}
\]
The corrector is accepted when the Euclidean norm of this augmented
residual is at most \(10^{-11}\), within at most \(20\) Newton
iterations, and the correction distance is at most \(2\Delta s\).
Each Newton update uses backtracking. Failed predictors cause
\(\Delta s\) to be halved; the minimum permitted step is
\(10^{-9}\). We start at \(\Delta s=10^{-4}\), increase the step
by a factor \(1.25\) after a corrector requiring at most three
iterations, and cap it at \(0.02\). After seven or more corrector
iterations the next step is multiplied by \(0.7\).

The maximum number of accepted steps is \(1500\). For this local
experiment we stop thirty accepted steps after the first change
from negative to nonnegative \(t\)-component of the tangent. This
condition was met at step \(625\), and the run stopped at step
\(655\). Neither the sampled minimum nor the tangent sign change
is an interval-certified turning point. All calculations use real
double precision. The saved full-tensor residuals were also checked
by a separate summation of the rank-one terms using extended
precision arithmetic.

The supplementary files contain the input coefficients, the
projection \(P\), the slice matrices \(Q,S\), the initial singular
values, all accepted-step diagnostics, and the factors at the
smallest recorded \(t\), the smallest full-tensor residual, and
the last accepted point. For every point we distinguish
\(\|P(f(\theta)-b(t))\|_2\),
\(\|f(\theta)-b(t)\|_2\), and
\(\rho(\theta)=\|f(\theta)-m\|_2/\sqrt{27}\).

\section{Completed epochs of the individual runs}
\label{app: stopping}

\Cref*{tab: stop rank} lists, for each run of the rank experiments (the $2\times2$ calibration and the $3\times3$ sweep), the number of epochs completed when the run was stopped.

\begin{table}[h!]
    \centering
    \small
    \caption{Rank experiments: epochs completed by each run at the interruption (runs $1$ to $7$). The programmed maximum is $10^{6}$ epochs.}
    \label{tab: stop rank}
    \begin{tabular}{c c c c c c c c c}
        \toprule
        Matrices & $r$ & $1$ & $2$ & $3$ & $4$ & $5$ & $6$ & $7$ \\
        \midrule
        $2\times2$ & $6$ & $263{,}340$ & $263{,}340$ & $263{,}340$ & $263{,}340$ & $263{,}340$ & $263{,}340$ & $263{,}340$ \\
         & $7$ & $304{,}690$ & $304{,}690$ & $304{,}690$ & $304{,}690$ & $304{,}690$ & $304{,}690$ & $304{,}690$ \\
        \addlinespace
        $3\times3$ & $16$ & $830{,}690$ & $830{,}690$ & $830{,}690$ & $830{,}690$ & $830{,}690$ & $830{,}690$ & $830{,}690$ \\
         & $17$ & $830{,}690$ & $830{,}690$ & $830{,}690$ & $830{,}690$ & $830{,}690$ & $830{,}690$ & $830{,}690$ \\
         & $18$ & $28{,}613$ & $28{,}662$ & $28{,}742$ & $28{,}566$ & $28{,}673$ & $28{,}739$ & $28{,}768$ \\
         & $19$ & $28{,}699$ & $28{,}573$ & $28{,}681$ & $28{,}712$ & $28{,}751$ & $28{,}607$ & $28{,}668$ \\
         & $20$ & $28{,}622$ & $28{,}757$ & $28{,}696$ & $28{,}541$ & $28{,}715$ & $28{,}604$ & $28{,}768$ \\
         & $21$ & $28{,}630$ & $28{,}509$ & $28{,}563$ & $28{,}576$ & $28{,}619$ & $28{,}612$ & $28{,}609$ \\
         & $22$ & $28{,}495$ & $28{,}727$ & $28{,}599$ & $28{,}544$ & $28{,}649$ & $28{,}608$ & $28{,}586$ \\
         & $23$ & $9{,}136$ & $9{,}146$ & $9{,}140$ & $9{,}138$ & $9{,}133$ & $9{,}145$ & $9{,}136$ \\
         & $24$ & $9{,}145$ & $7{,}285$ & $7{,}288$ & $7{,}280$ & $7{,}283$ & $7{,}282$ & $7{,}278$ \\
        \bottomrule
    \end{tabular}
\end{table}

\Cref*{tab: stop brank} lists the same information for the $126$ border-rank runs, together with the value of $\varepsilon$ reached at the interruption. The latter follows from the number of epochs: at epoch $e$ the schedule gives $\varepsilon=10^{-1}\cdot 10^{-2(e-1)/10^{6}}$.

\begin{table}[h!]
    \centering
    \small
    \setlength{\tabcolsep}{3.5pt}
    \caption{Border-rank experiments: epochs completed by each run at the interruption (runs $1$ to $7$), and range of $\varepsilon$ reached. The programmed maximum is $10^{6}$ epochs, at which $\varepsilon=10^{-3}$.}
    \label{tab: stop brank}
    \begin{tabular}{c c c c c c c c c c}
        \toprule
        $r$ & Exponents & $1$ & $2$ & $3$ & $4$ & $5$ & $6$ & $7$ & $\varepsilon$ at stop \\
        \midrule
        $16$ & $[-2,4]$ & $1{,}645$ & $1{,}647$ & $1{,}646$ & $1{,}648$ & $1{,}647$ & $1{,}647$ & $1{,}647$ & $0.0992$ \\
         & $[-3,5]$ & $1{,}647$ & $1{,}648$ & $1{,}648$ & $1{,}648$ & $1{,}646$ & $1{,}646$ & $1{,}647$ & $0.0992$ \\
         & $[-4,6]$ & $1{,}646$ & $1{,}647$ & $1{,}648$ & $1{,}646$ & $1{,}647$ & $1{,}648$ & $1{,}647$ & $0.0992$ \\
        \addlinespace
        $17$ & $[-2,4]$ & $1{,}628$ & $1{,}629$ & $1{,}627$ & $1{,}629$ & $1{,}628$ & $1{,}629$ & $1{,}629$ & $0.0993$ \\
         & $[-3,5]$ & $1{,}628$ & $1{,}628$ & $1{,}629$ & $1{,}628$ & $1{,}628$ & $1{,}628$ & $1{,}628$ & $0.0993$ \\
         & $[-4,6]$ & $1{,}628$ & $1{,}628$ & $1{,}629$ & $1{,}628$ & $1{,}629$ & $1{,}628$ & $1{,}628$ & $0.0993$ \\
        \addlinespace
        $18$ & $[-2,4]$ & $1{,}629$ & $1{,}629$ & $1{,}629$ & $1{,}628$ & $1{,}629$ & $1{,}629$ & $1{,}630$ & $0.0993$ \\
         & $[-3,5]$ & $1{,}628$ & $1{,}628$ & $1{,}629$ & $1{,}629$ & $1{,}630$ & $1{,}628$ & $1{,}629$ & $0.0993$ \\
         & $[-4,6]$ & $1{,}629$ & $1{,}629$ & $1{,}629$ & $1{,}629$ & $1{,}628$ & $1{,}628$ & $1{,}629$ & $0.0993$ \\
        \addlinespace
        $19$ & $[-2,4]$ & $1{,}550$ & $1{,}549$ & $1{,}550$ & $1{,}550$ & $1{,}550$ & $1{,}550$ & $1{,}549$ & $0.0993$ \\
         & $[-3,5]$ & $1{,}550$ & $1{,}549$ & $1{,}551$ & $1{,}549$ & $1{,}550$ & $1{,}550$ & $1{,}550$ & $0.0993$ \\
         & $[-4,6]$ & $1{,}549$ & $1{,}550$ & $1{,}549$ & $1{,}549$ & $1{,}547$ & $1{,}549$ & $1{,}264$ & $0.0993$ \\
        \addlinespace
        $20$ & $[-2,4]$ & $1{,}549$ & $1{,}550$ & $1{,}549$ & $1{,}551$ & $1{,}549$ & $1{,}549$ & $1{,}551$ & $0.0993$ \\
         & $[-3,5]$ & $1{,}549$ & $1{,}550$ & $1{,}549$ & $1{,}549$ & $1{,}550$ & $1{,}549$ & $1{,}549$ & $0.0993$ \\
         & $[-4,6]$ & $1{,}551$ & $1{,}548$ & $1{,}549$ & $1{,}548$ & $1{,}548$ & $1{,}549$ & $1{,}550$ & $0.0993$ \\
        \addlinespace
        $21$ & $[-2,4]$ & $3{,}835$ & $3{,}834$ & $3{,}835$ & $3{,}836$ & $3{,}834$ & $3{,}833$ & $3{,}834$ & $0.0983$ \\
         & $[-3,5]$ & $3{,}835$ & $3{,}835$ & $3{,}832$ & $3{,}834$ & $3{,}834$ & $3{,}835$ & $3{,}834$ & $0.0983$ \\
         & $[-4,6]$ & $3{,}834$ & $3{,}836$ & $3{,}834$ & $3{,}836$ & $3{,}835$ & $3{,}833$ & $3{,}835$ & $0.0983$ \\
        \bottomrule
    \end{tabular}
\end{table}

\end{document}